\documentclass{amsproc}

\newtheorem{theorem}{Theorem}[section]
\newtheorem{lemma}[theorem]{Lemma}
\newtheorem{proposition}[theorem]{Proposition}
\newtheorem{corollary}[theorem]{Corollary}

\theoremstyle{definition}
\newtheorem{definition}[theorem]{Definition}

\theoremstyle{remark}

\numberwithin{equation}{section}

\def\QED{{\nobreak\hfil\penalty50
  \hskip2em\hbox{}\nobreak\hfil
   \vrule height4pt width 4pt depth0pt 
   \parfillskip=0pt \finalhyphendemerits=0
      \par\medskip\goodbreak\noindent}}

  \newcommand{\cal}{\mathcal}

\input xy
\xyoption{all}

\usepackage{amsmath,url}
\usepackage{amscd}
\usepackage{amssymb}

\def\ra{\rightarrow}

\title{Quasi-categories vs Segal spaces}
\author{Andr\'e Joyal}
\address{Departement de Math\'ematiques, Universit\'e du Qu\'ebec \`a Montr\'eal} 
\email{joyal@math.uqam.ca}

\author{Myles Tierney}
\address{Departement de Math\'ematiques, Universit\'e du Qu\'ebec \`a Montr\'eal} 
\email{tierney@math.uqam.ca}

\thanks{The research of the  first author was supported by the NSERC}

\date{July 17th, 2006}

\begin{document}

\dedicatory{To Ross Street on the occasion of his 60th birthday}

\begin{abstract}
We show that  complete Segal spaces and Segal categories
are Quillen equivalent to quasi-categories. 
\end{abstract}

\maketitle

\section*{ Introduction}
Quasi-categories were introduced by 
Boardman and Vogt  in their work on homotopy invariant
algebraic structures \cite{BV}. 
They are often called weak Kan complexes in the literature.
The category  of simplicial sets ${\bf S}$ 
admits a Quillen model structure in which the fibrant objects are
the quasi-categories by a result of the first author in \cite{J2}. 
We call this model structure, the {\it model structure for quasi-categories}. 
The theory of quasi-categories has applications to homotopical
algebra and higher category theory, see \cite{J3}, \cite{Lu1} and \cite{Lu2}.
Complete Segal spaces were introduced by Rezk
in his work on the homotopy theory of homotopy theories  \cite{Rez}.
The category of bisimplicial sets ${\bf S}^{(2)}$ 
admits a Quillen model structure in which the fibrant objects
are the complete Segal spaces.
We call this model structure the {\it model structure for complete Segal spaces}. 
The main result of this paper is to establish a Quillen equivalence
between the model category for quasi-categories
and that for complete Segal spaces:
$$p_1^*:{\bf S} \leftrightarrow {\bf S}^{(2)}:i_1^*.$$
The functor $i_1^*$ associates to a bisimplicial set $X$
its first row $X_{\star 0}$. The theorem implies 
that the first row of a complete Segal space
contains all the homotopy information
about the space. We also describe a Quillen equivalence
in the opposite direction, 
 $$t_!:{\bf S}^{(2)} \leftrightarrow {\bf S}:t^!,$$
where the functor $t_!$ associates
to a bisimplicial set $X$ a total simplicial set $t_!(X)$.
Segal categories were first introduced by Dwyer, Kan
and Smith \cite{DKS}, where they are called special $\Delta^o$-diagrams
of simplicial sets .
The theory of Segal categories was extensively developed by 
Hirschowitz and Simpson for application to algebraic geometry.
They show that the category of precategories admits a model structure
in which the fibrant objects are Segal categories.
We call this model structure the {\it model structure for Segal
categories}. 
By a theorem of Bergner  in \cite{B2}, the inclusion functor $\pi^*:{\bf PCat}\subset  {\bf S}^{(2)}$
induces a Quillen equivalence between
the model structure for Segal
categories and the model structure for 
complete Segal spaces,
 $$\pi^*:{\bf PCat} \leftrightarrow {\bf S}^{(2)}:\pi_*.$$
By combining these equivalences,
we obtain an equivalence
between the model category for quasi-categories
and that for Segal categories: 
$$q^*:{\bf S} \leftrightarrow {\bf PCat}:j^*.$$
The functor $j^*$ associates to a precategory $X$
its first row $X_{\star 0}$. The theorem implies 
that the first row of a fibrant Segal category
contains all the homotopy information
about the Segal category. 
We also describe a Quillen equivalence
in the opposite direction, 
 $$d^*:{\bf PCat} \leftrightarrow {\bf S}:d_*,$$
where the functor $d^*$ associates
to a precategory $X$ its diagonal simplicial set.

 \medskip

The paper has five sections, one addendum, one appendix and one epilogue.
 In the first section we give a brief description of the relevant aspects
of the model structure for  quasi-categories. 
 The fibrant objects are the quasi-categories and the acyclic maps are called
weak categorical equivalences. 
In the second section we introduce two Reedy model structures
on the category of bisimplicial sets, respectively called
the vertical and the horizontal model structures.
The acyclic maps are the column-wise weak homotopy
equivalences in the vertical structure 
but they are the row-wise weak categorical equivalences
in the horizontal.
The interaction between the two model structures 
is one of the main technical tools of the paper. We define a total space functor
from bisimplicial sets to simplicial sets
and show that it is a left Quillen functor with values in the model category for quasi-categories
for both the vertical and the horizontal model structures.
In the third section we show that 
every row of a Segal space is a quasi-category.
In section four we show that the functor which associates to a bisimplicial set its first row
induces a Quillen equivalence
between complete Segal spaces
and quasi-categories. 
We also show that the total space functor induces
a Quillen equivalence in the opposite direction.
In section five we show that the functor which associates to a precategory 
its first row induces a Quillen equivalence
between Segal categories
and quasi-categories. 
We also prove that the functor which associates to a precategory its diagonal
induces a Quillen equivalence in the opposite direction.
In the addendum we show that 
the model category for complete Segal
spaces can be obtained from the model category
for quasi-categories by a method developped by 
Rezk, Schwede and Shipley
in their paper on  model categories and functors \cite{RSS}.
In the epilogue, we discuss other models of homotopy theories
and other Quillen equivalences.

\tableofcontents

\section {The model structure for quasi-categories}

\medskip

The results of this section are taken from \cite{J1} and \cite{J2}. 
Quasi-categories were introduced by 
Boardman and Vogt  in their work on homotopy invariant
algebraic structures \cite{BV}. 
The category of simplicial sets ${\bf S}$ admits a model structure in which
 the fibrant objects are the quasi-categories \cite{J2}. 
 We call it the {\it model structure for quasi-categories}.

\bigskip

See the appendix for notation, and basic results.
We denote by ${\bf Set}$ the category of sets,  by ${\bf S}$
the category of simplicial sets 
and by ${\bf Cat}$ the category of small categories.

\medskip

Recall that
an arrow $u:A\to B$ in a category is said to have the
{\it left lifting property} with respect to another arrow $f:X\to Y$, or that
$f$ has the {\it right lifting property} with respect to $u$,
if every
commutative square
$$\xymatrix{
A \ar[d]_u \ar[r]^x &  X  \ar[d]^f  \\
B \ar[r]_y  \ar@{-->}[ur]       &Y    }$$
has a diagonal filler $d:B\to X$ (that is, $du=x$ and $fd=y$).
We shall denote this relation by 
$u \pitchfork f$.

\medskip

\medskip

We first
recall the classical model 
structure on the category ${\bf S}$.
The notion of weak homotopy equivalence
between simplicial sets is usually defined by using the
geometric realisation functor $ {\bf S} \to {\bf Top}$.
An alternative definition uses Kan complexes and the 
homotopy category ${\bf S}^{\pi_0}$ introduced
by Gabriel and Zisman \cite{GZ}.
Recall that a simplicial set $X$ is called a {\it Kan complex}
if every  horn
$\Lambda^k[n] \rightarrow X$ has a filler $\Delta[n]\rightarrow X$.
If $A,B\in {\bf S}$ let us put
$$\pi_0(A,B)=\pi_0(B^A).$$
If we apply  the functor $\pi_0$
to  the composition map $C^B\times B^A\to C^A$
we obtain a
composition law
$${\pi_0}(B,C)\times {\pi_0}(A,B)\to  {\pi_0}(A,C)$$
for a category ${\bf S}^{\pi_0}$,
where ${\bf S}^{\pi_0}(A,B)=\pi_0(A,B)$.
A map of simplicial sets  is called
a {\it homotopy equivalence} if it is invertible
in the category ${\bf S}^{\pi_0}$.
A map $u:A\to B$ is called
a {\it weak homotopy equivalence} 
the map 
$$\pi_0(u,X):\pi_0(B,X)\to \pi_0(A,X)$$
is bijective
for every Kan complex $X$.
This notion of weak equivalence is equivalent to the usual notion
which is defined via the geometric realisation functor.

\medskip

A map of simplicial sets is called a {\it Kan fibration}
if it has the right lifting property with respect to every horn inclusion 
$\Lambda^k[n] \subset \Delta[n]$. 
The following theorem describes the {\it classical model structure} on
the category $\bf S$.

\begin{theorem} [Quillen, see \cite{Q}]     \label{classicalmodel}
The category of simplicial sets $\bf S$ admits a model structure $({\cal C}_0,{\cal W}_0,{\cal F}_0)$
in which a cofibration is a monomorphism,
a weak equivalence is a weak homotopy equivalence
and a fibration is a Kan fibration.
The model structure is cartesian closed and proper.
\end{theorem}

See Definition \ref{defmodelcat} in the appendix for the notion of model structure.
See \cite{JT2} for a purely combinatorial proof of the theorem. 
We call a map of simplicial sets a {\it trivial fibration}
if it has the right lifting property with respect to every 
monomorphism. We note that this notion can be defined
in any topos, see Definition \ref{trivfib}. 
The acyclic fibrations of the classical model structure are the trivial fibrations.
For every $n\geq 0$, we denote by $\delta_n$ the map $\partial \Delta[n] \to \Delta[n]$
defined by the inclusion $\partial \Delta[n] \subset \Delta[n]$.

\begin{proposition} \label{saturatedmonosset} \cite{GZ}
The class of monomorphisms in the category $\bf S$
is generated as a saturated class by
the set of inclusions 
$$\delta_n: \partial \Delta[n] \subset \Delta[n], \quad {\rm for}\enspace n\geq 0.$$
A map of simplicial sets is a trivial fibration
iff it has the right lifting property
with respect to $\delta_n$ for
every $n\geq 0$.
\end{proposition}

For the notion of a saturated class, see \ref{defsaturated}.

\begin{definition}\label{anodyne def} \cite{GZ}
A map of simplicial sets is said to be {\it anodyne} 
if it belongs to the saturated class generated by the horns  $\Lambda^k[n] \subset \Delta[n]$
for $n\geq 1$ and $0\leq k\leq n$.
\end{definition}

A map
of simplicial set is anodyne iff it is an acyclic cofibration.

\medskip

Before describing the model structure for quasi-categories it is good to 
describe a related  model 
structure on ${\bf Cat}$. 
We call a functor $p:X\to Y$
a {\it quasi-fibration}  if
for every object $a\in X$ and every isomorphism $g\in Y$
with source $p(a)$ there exists an isomorphism $f\in X$
with source $a$ such that $p(f)=g$.
A functor $p:X\to Y$
is a quasi-fibration iff it has the right lifting property
with respect to the inclusion $\{0\}\subset J$,
where $J$ is the groupoid generated by one isomorphism $0\rightarrow 1$.
We say that a functor $A\to B$ is {\it monic on objects} (resp.  {\it surjective on objects})
if the induced map $Ob (A)\to Ob(B)$ is monic  (resp. surjective).

\begin{theorem} \label{Catmodel}  \cite{JT1} \cite{Rez} The category ${\bf Cat}$
admits a model structure 
in which a cofibration is a functor monic on objects,
a weak equivalence is an equivalence of categories
and a fibration is a quasi-fibration.
The acyclic fibrations are the equivalences surjective on objects.
The model structure is proper and cartesian closed.
Every object is fibrant and cofibrant.
\end{theorem}

We call this model structure the {\it natural model structure}  on ${\bf Cat}$.
The notion of cartesian closed model category is defined
in \ref{defcartclosedmod}.

\bigskip

The category $\Delta $ is a full subcategory of $ {\bf Cat}$.
The {\it nerve}  of a category $C\in  {\bf Cat}$ is 
the simplicial set $NC$ obtained by putting $(NC)_n={\bf Cat}([n],C)$
for every $n\geq 0$.
The nerve functor
$N: {\bf Cat}\rightarrow  {\bf S}$ is full and faithful and we shall regard it as an 
inclusion $N: {\bf Cat}\subset {\bf S}$
by adopting the same notation for a category and its nerve.
The functor $N$
has a left adjoint 
$$\tau_1:{\bf S} \to  {\bf Cat}.$$
We say that
$\tau_1 X$ is the {\it fundamental category} of a simplicial set $X$.
The  fundamental groupoid $\pi_1 X$ is obtained 
by inverting the arrows of $\tau_1 X$.

\bigskip

We shall say that a
horn
$\Lambda^k[n] \subset \Delta[n]$ is {\em inner}
if $0 < k < n$. 
The following definition is due to Boardman and Vogt.

\begin{definition} \label{defquasicat} \cite{BV} A simplicial set $X$ 
is called a {\it quasi-category} if 
every inner horn $\Lambda^k[n] \to X$ has a filler $\Delta[n] \to X$. 
A {\it map} of quasi-categories
is a map of simplicial sets.
\end{definition}

Quasi-categories are often called {\it weak Kan complexes} in the literature.
The nerve of a category and a Kan complex are examples.
We shall denote by ${\bf QCat}$
the category of quasi-categories;
it is a full subcategory of $\bf S$.

\medskip

The next step is to introduce an appropriate notion
of equivalence for quasi-categories. If $A$ is a simplicial set,
we shall denote by $\tau_0 A$ 
the set of isomorphism classes of objects
of the fundamental category $\tau_1 A$.
The functor $\tau_1: {\bf S} \to {\bf Cat}$ preserves finite products \cite{GZ},
hence also the functor $\tau_0: {\bf S} \to  {\bf Set}.$
If $A,B\in {\bf S}$ let us put
$$\tau_0(A,B)=\tau_0(B^A).$$
If we apply the functor $\tau_0$
to the composition map $C^B\times B^A\to C^A$
we obtain a
composition law
$${\tau_0}(B,C)\times {\tau_0}(A,B)\to  {\tau_0}(A,C)$$
for a category ${\bf S}^{\tau_0}$,
where ${\bf S}^{\tau_0}(A,B)=\tau_0(A,B)$.
We shall say that a map of simplicial sets 
is a {\it categorical equivalence} 
if it is invertible
in the category ${\bf S}^{\tau_0}$.
If $X$ and $Y$ are quasi-categories,
a categorical equivalence $X\to Y$
is called an {\it equivalence of quasi-categories}. 
We shall say that a map of simplicial sets $u:A\to B$
is a {\it weak categorical equivalence} 
if the map 
$$\tau_0(u,X):\tau_0(B,X)\to \tau_0(A,X)$$
is bijective
for every quasi-category $X$.

\medskip

The notion of equivalence between quasi-categories
has another description which we don' t need
but which is good to know. If $X$ be a simplicial set, we shall denote by $X(a,b)$
the fiber at $(a,b)\in X_0\times X_0$
of the projection $(p_0,p_1):X^I\to X\times X$  
defined from the inclusion $\{0, 1\}\subset I$.
The simplicial set $X(a,b)$
is a Kan complex if  $X$ is a quasi-category.

\begin{definition}\label{ffmapofqcat} \cite{J2}
We shall say that a map between
simplicial sets $u:A\to B$ is {\it essentially surjective}
if the map $\tau_0(u):\tau_0 A\to \tau_0 B$
is surjective.
We shall say that a map between quasi-categories
$f:X\to Y$ is {\it fully faithful} if 
the map
$$X(a,b)\to Y(fa,fb)$$
induced by $f$ is a weak homotopy equivalence for every pair $a,b\in X_0$.
\end{definition}

\begin{theorem}\label{equivisfullyfaithandesssur}  \cite{J2} A map  between quasi-categories 
is an equivalence iff it is fully faithful and essentially surjective.
\end{theorem}

\begin{definition}  \label{defquasi-fib}  \cite{J2} We shall say that a map of simplicial sets  
is a {\it quasi-fibration} if it has the right
lifting property with respect to every monic weak categorical equivalence.
\end{definition}

\begin{theorem} \label{QCatmodel}\cite{J2}
The category of simplicial sets $\bf S$ admits a model structure \break
$({\cal C}_1,{\cal W}_1,{\cal F}_1)$ in which a cofibration is a 
monomorphism, a weak equivalence is a 
weak categorical equivalence
and a fibration is a quasi-fibration.
The acyclic fibrations are the trivial fibrations.
The fibrant objects are the quasi-categories.
The model structure is cartesian closed and left proper.
\end{theorem}

We shall say that it is the {\it model structure for quasi-categories}.
The notion of cartesian closed model category is defined in \ref{defcartclosedmod}.
The quasi-fibrations between quasi-categories have
a simple description. It is based on 
the following notion:

\begin{definition} \cite{J2} \label{mid cofibration def} We shall say that a map of simplicial sets is 
{\it mid anodyne} 
if it belongs to the saturated class generated by the inner horns  $\Lambda^k[n] \subset \Delta[n]$.
We shall say that a map is a {\it mid fibration} 
if it has the right
lifting property with respect to every  inner horn $\Lambda^k[n] \subset \Delta[n]$.
\end{definition}

For the notion of a saturated class, see \ref{defsaturated}.

\begin{proposition}\label{midcofibration}  \cite{J2} Every mid anodyne map is a weak categorical equivalence bijective on vertices. The functor $\tau_1:{\bf S}\to {\bf Cat}$ 
takes a mid anodyne map
to an isomorphism of categories.
\end{proposition}

\begin{proposition} \label{midcofibrationfibration} \cite{J2} If  ${\cal B}$ is the class of mid fibrations
and ${\cal A}$ is the class of mid anodyne maps
then  the pair $({\cal A}, {\cal B})$
is a weak factorisation system in the category  ${\bf S}$
\end{proposition} 

For the notion of weak factorisation system, see \ref{defweakfact}.

\medskip

Let us regard the groupoid $J$ as a simplicial set via the nerve functor.

\begin{proposition}  \label{qfbetqc} \cite{J2} Every quasi-fibration is a mid fibration.
Conversely, a mid fibration between quasi-categories $p:X\to Y$
is a quasi-fibration iff the following equivalent conditions are satisfied:
\begin{itemize}
\item{} $p$ has the right lifting property
with respect to the inclusion $\{0\}\subset J$
\item{} the functor $\tau_1 p:\tau_1X\to \tau_1Y$
is a quasi-fibration. 
\end{itemize}
\end{proposition}

\medskip

\begin{proposition}\label{Qpair1}\cite{J2} The pair of adjoint
functors $$\tau_1:{\bf S}\leftrightarrow  {\bf Cat}:N$$
is a Quillen pair between the model structure for quasi-categories
and  the natural model structure on ${\bf Cat}$.
A functor $u:A\to B$ in $\bf Cat$
is a quasi-fibration (resp. an equivalence) 
iff the map $Nu:NA\to NB$ is a quasi-fibration (a weak categorical
equivalence) in $\bf S$.
\end{proposition}

It follows that the functor  $\tau_1:{\bf S}\to {\bf Cat}$ 
takes a weak categorical equivalence
to an equivalence of categories.

\medskip

\begin{proposition}\label{Qpair2} \cite{J2} The classical model
structure on ${\bf S}$
is a Bousfield localisation of the model structure for
quasi-categories. 
 \end{proposition}

For the notion of Bousfield localisation, see Definition \ref{DefBousfieldLocalisationA}. 
Thus, every weak categorical equivalence
is a weak homotopy equivalence
and every Kan fibration 
is a quasi-fibration. Conversely,
a map between Kan complexes is a weak homotopy equivalence (resp, a Kan fibration)
iff it is a weak categorical equivalence (resp. a quasi-fibration).

\medskip

Let ${\bf Kan }$ be the category of Kan complexes.

\begin{proposition}  \label{Jfunctoradjoint}  \cite{J2}
The inclusion functor ${\bf Kan }\subset  {\bf QCat}$
has
a right adjoint $$J:{\bf QCat}\to {\bf Kan }.$$
The Kan complex $J(X)$
is the largest sub-Kan complex of a quasi-category $X$.
The functor $J$  takes an equivalence of quasi-categories
to a homotopy equivalence and a quasi-fibration
to a Kan fibration.
\end{proposition}

If $X$ is a simplicial set, then the contravariant functor
$A\mapsto X^A$ is right adjoint to itself.
If $X$ is a quasi-category,  then so is the simplicial set $X^A$.
We shall denote by $X^{(A)}$ be
the full simplicial subset of $X^A$
whose vertices are the maps
$f:A\to X$ such that $f(A)\subseteq J(X)$.
The contravariant functor $A\mapsto J(X^A)$ 
is a subfunctor of the contravariant functor $A\mapsto X^A$.

\begin{proposition}\label{Jdual2}  \cite{J2} Let $X$ be a quasi-category.
Then the contravariant functors $A\mapsto J(X^A)$ and $A\mapsto X^{(A)}$
are mutually right adjoint.
The contravariant functor $A\mapsto X^{(A)}$
takes a weak homotopy equivalence to an equivalence 
of quasi-categories.
\end{proposition}

\medskip

Consider the functor $k:\Delta\to {\bf S}$.
defined by putting $k[n]=\Delta'[n]$
for every $n\geq 0$, where $\Delta'[n]$ denotes
the (nerve of the) groupoid freely generated by
the category $[n]$. We denote by 
$k^!: {\bf S} \to {\bf S}$ the functor defined by putting 
$$k^!(X)_{n}={\bf S}( \Delta'[n],X)$$
for every simplicial set $X$ and for every $n\geq 0$. The functor $k^!$ 
has a left adjoint $k_!: {\bf S}\to {\bf S}$
which is the left Kan extension of the functor
$k:\Delta\to {\bf S}$ along the Yoneda functor $y:\Delta\to {\bf S}$.

\begin{proposition}\label{kandpi} For every $X\in {\bf S}$,
we have $\tau_1k_!X=\pi_1X$.
\end{proposition}

\noindent
{\bf Proof}:  The result is obvious if $X=\Delta[n]$.
The general case follows, since the functors 
$X\mapsto \tau_1k_!X$ and $X\mapsto \pi_1X$
are cocontinuous. \QED

\begin{theorem}\label{Qpair3} \cite{J2} The pair of adjoint
functors $$k_!:{\bf S}\leftrightarrow  {\bf S}:k^!$$
is Quillen  pair
$ ({\cal C}_0, {\cal W}_0,{\cal F}_0) \leftrightarrow  ({\cal C}_1, {\cal W}_1,{\cal F}_1).$
\end{theorem}

The pair $(k_!,k^!)$
is actually a homotopy colocalisation in the sense of \ref{Localisation1}.

\medskip

The inclusion $\Delta[n]\subseteq \Delta'[n]$
is natural in $[n]\in \Delta$.
Hence it defines a natural
transformation $y\to k$, where $y$ is the Yoneda functor $\Delta\to {\bf S}$.
The transformation has a unique extension
$Id\to k_!$, where $Id$ is the identity functor of ${\bf S}$.
There is a corresponding adjoint transformation 
$k^!\to Id$. It is easy to verify that
the map
$X\to k_!(X)$ is monic and a weak homotopy equivalence
for every $X\in {\bf S}$. Dually, 
the map
$k^!(X)\to X$ is a trivial fibration
for every Kan complex $X$. 
If $X$ is a quasi-category, then the map $k^!(X)\to X$ induces a map $k^!(X)\to J(X)$.

\begin{proposition}  \label{kfunctorJfunctor} \cite{J2}
The natural map $k^!(X)\to J(X)$
is a trivial fibration for every quasi-category $X$.
\end{proposition}

Let $X$ be a quasi-category. 
Then the canonical map $X^{k_!(A)}\to X^{(A)}$
defined from  the inclusion $A\subseteq k_!(A)$
is a categorical equivalence for every simplicial set $A$.

\bigskip

The following combinatorial result will be used in section 4. 

\medskip

For $n > 0$, the {\em $n$-chain} $I_n\subseteq \Delta[n]$ is defined to be the
union  of the 
edges $(i, i+1) \subseteq \Delta[n]$ for $0 \leq i \leq n-1$. 
Let us put $I_0=1$.

\begin{lemma}\label{midado2}\cite{J2} The inclusion $I_n\subseteq \Delta[n] $
is mid anodyne .\end{lemma}

\bigskip

The following notion will be used in section 5. 

\medskip

Recall that a functor $u:A\to B$ is said to be {\it conservative}
if the implication
$$u(f)\quad  {\rm invertible} \quad \Rightarrow \quad
 f \quad  {\rm invertible}$$
is true for every arrow $f\in A$.

\begin{definition} \cite{J2} \label{conservativedef} We shall say that a map of simplicial sets $u:A\to B$
is {\it conservative} if the functor
$\tau_1(u):\tau_1 A\to \tau_1 B$ is conservative. 
\end{definition}

\begin{theorem} \label{conservativeprojection} \cite{J2}
If $X$ is a quasi-category and $A$ is a simplicial set,
then 
the projection 
$$X^A\to X^{A_0}$$
is conservative.
\end{theorem}

\begin{proposition} \label{conservativefibKanfibration}\cite{J2} Let $p:X\to Y$ be a 
conservative quasi-fibration
between quasi-categories. 
If $Y$ is a Kan complex, then $X$ is a Kan complex and $p$ is a Kan fibration.
\end{proposition}

\begin{corollary} \label{basechangeconservativeJ} \cite{J2} The base change in
${\bf QCat}$
of a conservative quasi-fibration between quasi-categories  is conservative.
\end{corollary}

\medskip

\section {The vertical and horizontal Reedy model structures}

In this section we introduce two model
structures which play an important role
in the paper.  They are called the vertical and the horizontal model structures.
Each is a Reedy model structure 
on the category of bisimplicial sets 
associated to a model structure on the category of simplicial sets.
In the vertical  model structure the acyclic maps 
are the column-wise weak homotopy equivalences
but  in the horizontal they
are the row-wise weak categorical equivalences.
The vertical model structure is associated to  the classical model structure
on the the category of simplicial sets while
the horizontal is associated to the quasi-category model structure.
We define a ``total space" functor $t_!:{\bf S}^{(2)} \to {\bf S}$ 
and show that it
is a left Quillen functor with respect to both
model structures.

\bigskip

A {\it bisimplicial set} is a contravariant functor  
$\Delta \times  \Delta \to {\bf Set}$.
We shall denote the category
of bisimplicial sets by ${\bf S}^{(2)}$.
A {\it simplicial space} is a contravariant functor  
$\Delta \to {\bf S}$. 
We regard a simplicial space $X$
as a bisimplicial set
by putting $X_{mn}=(X_m)_n$ for every $m,n\geq 0$.
Conversely, we regard a bisimplicial set $X$
as a simplicial space by putting $X_m=X_{m\star}$ for every $m\geq 0$.
The {\em box product} $A \Box B$ of two simplicial sets 
$A$ and $B$ is the bisimplicial set obtained 
by putting 
$$(A \Box B)_{mn} = A_m \times B_n$$
for every $m,n\geq 0$.  This defines a functor of 
two variables $\Box :{\bf S}\times {\bf S} \to {\bf S}^{(2)}$.
The functor is divisible on both sides. 
This means that the
functor $A\Box(-):  {\bf S}   \to  {\bf S}^{(2)}$
admits a right adjoint $A\backslash (-): {\bf S}^{(2)} \to {\bf S}$
for every simplicial set $A$.
If $X\in {\bf S}^{(2)}$,  then a simplex $\Delta[n]\to A\backslash X$
is a map $A\Box \Delta[n]\to X$.
The simplicial set
$\Delta[m] \backslash X$ is the $m$th column $X_{m*}$ of $X$.
Dually, the functor $(-)\Box B:  {\bf S}   \to  {\bf S}^{(2)}$
admits a right adjoint $(-)/B:  {\bf S}^{(2)} \to  {\bf S} $
for every simplicial set $B$.
If $X\in {\bf S}^{(2)}$,  then a  simplex $\Delta[m]\to X/B$
is a map $\Delta[m] \Box B \to X$.
The simplicial set
$X\slash\Delta[n] $ is the $n$th row $X_{*n}$ of  $X$.
If $X\in {\bf S}^{(2)}$ and $A,B\in  {\bf S}$, 
there is a bijection between
the following three kinds of maps
$$A \Box B\to X, \quad   \quad B\to A \backslash X,\quad    \quad A\to X\slash B.$$
Hence the contravariant functors $A\mapsto A \backslash X$
and $B\mapsto B \backslash X$ are mutually right adjoint.

\bigskip

 If $u:A\to B$ 
and $v:S\to T$ are maps of simplicial sets
we shall denote by $u   \Box' v$
the map 
$$A \Box  T  \sqcup_{A  \Box   S}  B \Box  S  \longrightarrow B  \Box    T
$$
obtained from the commutative square 
$$
\xymatrix{
A \Box S \ar[d] \ar[r]&  B  \Box  S  \ar[d]\\
A  \Box  T \ar[r]    &B  \Box    T .   }
$$
This defines a functor of two variables
$$\Box' :{\bf S}^I\times {\bf S}^I\to ({\bf S}^{(2)})^I,$$
where ${\cal E}^I$ denotes the category of arrows
of a category ${\cal E}$. 
The functor $\Box' $ is divisible on both sides.
If $u:A\to B$ is map in ${\bf S}$
and $f:X\to Y$ is a map in ${\bf S}^{(2)}$,
we denote by $\langle u\backslash\ f\rangle$
the map 
$$B\backslash X \to B\backslash Y  \times_{ A\backslash Y }  A\backslash X$$
obtained from the
square
$$\xymatrix{
B\backslash X \ar[r]\ar[d]&  A\backslash X \ar[d] \\
B\backslash Y   \ar[r]&  A\backslash Y .
 }
 $$
The functor  $f\mapsto \langle u\backslash f\rangle$
is right adjoint to the functor $v\mapsto u\Box' v$.
Dually, if $v:S\to T$ is map in ${\bf S}$,
we shall denote by $\langle f/v\rangle$
the map 
$$X/T \to Y/T \times_{ Y/S}  X/S$$
obtained from the
commutative square
$$\xymatrix{
X/T \ar[r]\ar[d]&  X/S\ar[d] \\
Y/T  \ar[r]&  Y/S .
 }
$$
The functor  $f\mapsto \langle f\backslash v\rangle$
is right adjoint to the functor $u\mapsto u\Box' v$.

\begin{proposition}  \label{lift2} 
For any triple of maps $u\in {\bf S}$, $v \in {\bf S}$ and $f\in {\bf S}^{(2)}$ we have
$$(u  \Box' v) \pitchfork f \quad \Longleftrightarrow \quad u\pitchfork \langle f/v \rangle \quad 
\Longleftrightarrow  \quad v\pitchfork \langle u\backslash f \rangle.$$
\end{proposition}

This follows from \ref{lift2A}.

\bigskip

Let us denote by $\delta_n$ the inclusion $\partial \Delta[n] \subset \Delta[n]$.
Recall that a map of bisimplicial sets
is said to be a {\it trivial fibration} if it has the right
lifting property with respect to every monorphism.

\begin{proposition} \label{saturatedmonobisset} \cite{JT2} The
class  of monomorphisms in the category ${\bf S}^{(2)}$ is generated as a saturated class
by the set of inclusions
$$\delta_m\Box' \delta_n: (\partial  \Delta[m] \Box \Delta[n]) \cup  (\Delta[m] \Box \partial \Delta[n]) \subset \Delta[m] \Box \Delta[n],\quad {\rm for} \quad  m,n\geq 0.$$
A map of bisimplicial sets is 
a trivial fibration
iff it has the right lifting property
with respect to the map $\delta_m\Box' \delta_n$ for every $ m,n\geq 0$.
\end{proposition}

\begin{proposition} \label{trivfibbisimpset} 
A map of bisimplicial sets $f:X\to Y$ is a trivial fibration 
iff the following equivalent conditions are satisfied:
\begin{itemize}
\item the map $\langle \delta_m\backslash f\rangle$ is a trivial fibration
for every $m\geq 0$;
\item  the map $\langle u\backslash f\rangle$ is a trivial fibration
for every monomorphism $u\in {\bf S}$; 
\item the map $\langle f/ \delta_n\rangle$ is a trivial fibration
for every $n\geq 0$;
\item  the map $\langle f/v \rangle$ is a trivial fibration
for every monomorphism $v\in {\bf S}$.
\end{itemize}
\end{proposition}

\noindent{\bf Proof}: Let us show that $f$ is a trivial fibration iff the map
$\langle \delta_m\backslash f\rangle$ is a trivial fibration
for every $m\geq 0$. But the map $\langle \delta_m\backslash f\rangle$ is a trivial fibration
iff we have $\delta_n\pitchfork \langle \delta_m\backslash f\rangle$
for every $n\geq 0$ by \ref{saturatedmonosset}.  But the condition 
$\delta_n\pitchfork \langle \delta_m\backslash f\rangle$
is equivalent to the condition $(\delta_m\Box' \delta_n)\pitchfork f$
by \ref{lift2}. Hence the result follows from \ref{saturatedmonobisset}.
Let us now show that first condition implies the fourth. 
But the map $\langle f/v \rangle$ is a trivial fibration
iff we have $\delta_m\pitchfork \langle f/v \rangle$
for every $m\geq 0$ by \ref{saturatedmonosset}. 
But the condition 
$\delta_m\pitchfork \langle f/v \rangle$
is equivalent to the condition $v\pitchfork \langle \delta_m\backslash f\rangle$
by \ref{lift2}. This proves the result since $\langle \delta_m\backslash f\rangle$
is a trivial fibration by hypothesis and since $v$ is monic.
The rest of the equivalences are proved similarly.
\QED

We shall use the following simplicial enrichement 
of the category of bisimplicial sets ${\bf S}^{(2)}$.
The functor  $i_2:\Delta \to  \Delta\times \Delta$ 
defined by putting $i_2([n])=([0],[n])$ 
is right adjoint to the second projection 
$p_2:\Delta \times \Delta \to \Delta $.
Hence the functor $i_2^*:{\bf S}^{(2)} \to {\bf S}$ 
is right adjoint to the functor $p_2^*$.
If $X$ is a bisimplicial set, then $i_2^*(X)$
is the first column of $X$.
If $A$ is a simplicial set, then $p_2^*(A)=1\Box A$.
For any pair of bisimplicial sets $X$ and $Y$,
let us put 
$$ Hom_2(X,Y) =i_2^*(Y^X).$$ 
This defines an enrichment
of the category ${\bf S}^{(2)}$
over the category ${\bf S}$.

\begin{proposition} \label{bisimpisenriched} The enriched category $({\bf S}^{(2)},Hom_2)$ 
admits tensor and cotensor products. The tensor product of a simplicial space $X$
by a simplicial set $A$ is the simplicial space $X\times_2 A=X\times p_2^*(A)$
and the cotensor product is the simplicial space $X^{p_2^*(A)}=X^{1\Box A}.$
\end{proposition}

\medskip

We say that a map of bisimplicial sets
$f:X\to Y$
is a {\it column-wise weak homotopy equivalence}
 if the map $\Delta[m]\backslash f=f_{m}:X_{m}\to Y_{m}$ 
is a weak homotopy equivalence for every $m\geq 0$.
We say that 
$f:X\to Y$  is a {\it vertical fibration}, or a {\it v-fibration},  if 
 the map
$\langle \delta_m\backslash f \rangle$
is a Kan fibration for every $m\geq 0$
(a notion of horizontal fibration
will be considered later).
We say that a bisimplicial set $X$ 
is {\it v-fibrant} if  the map $X\to 1$
is a v-fibration.

\medskip

Let $h^k_n$ denotes the inclusion $\Lambda^k[n] \subset \Delta[n] $.

\begin{proposition} \label{Reedyfibrationvertically} A map of bisimplicial sets $f:X\to Y$ is 
a v-fibration
iff it satisfies the following equivalent conditions:
\begin{itemize}
\item  the map
$\langle \delta_m  \backslash f \rangle $
is a Kan fibration for every $m\geq 0$;
\item  the map $\langle u  \backslash f \rangle $
is a Kan fibration for every monomorphism $u$;
\item  the map
$\langle f/ h^k_n \rangle $
is a trivial fibration for every $n>0$ and $0\leq k\leq n$;
\item  the map $\langle f/ v \rangle $
is a trivial fibration for every anodyne map $v\in {\cal S}$.
\end{itemize}
\end{proposition}

\noindent
{\bf Proof}: The equivalences (i)$\Leftrightarrow$(ii)$\Leftrightarrow$(iv)
follow from \ref{Reedyfib1}. The implication (iv)$\Rightarrow$(iii) is obvious.
Let us prove the implication (iii)$\Rightarrow$(ii).
Let us show that the map $\langle \delta_m  \backslash f \rangle $
is a Kan fibration for every $m\geq 0$.
For this, it suffices to show that we have $h^k_n \pitchfork \langle \delta_m  \backslash f \rangle$
for every $n>0$ and $0\leq k\leq n$.
But the condition $h^k_n\pitchfork \langle \delta_m  \backslash f \rangle$
is equivalent to the condition $\delta_m \pitchfork \langle  f /h^k_n\rangle$
by \ref{lift2}. But we have $\delta_m \pitchfork \langle  f /h^k_n\rangle$
since $\langle f/ h^k_n \rangle $
is a trivial fibration by assumption.
The implication (iii)$\Rightarrow$(ii) is proved.
\QED

\medskip

The following theorem describes the {\it vertical model structure} on ${\bf S}^{(2)}$.

\medskip

\begin{theorem}  \label{verticalReedymodel} \cite{Ree}
The category $({\bf S}^{(2)}, Hom_2)$ 
admits a simplicial model structure $({\cal C}^v_0, {\cal W}^v_0,{\cal F}^v_0)$ 
in which the cofibrations are the monomorphisms,
the weak equivalence are the column-wise weak homotopy equivalences
and the fibrations are the vertical fibrations.
The model structure is proper and cartesian closed.
The acyclic fibrations are the trivial  fibrations.
\end{theorem}

\noindent{\bf Proof}:  The model structure $({\cal C}^v_0, {\cal W}^v_0,{\cal F}^v_0)$ 
is the Reedy model structure $({\cal C}'_0, {\cal W}'_0,{\cal F}_0)$ associated to  
the classical model structure $({\cal C}_0, {\cal W}_0,{\cal F}_0)$ of  \ref{classicalmodel}  on 
${\bf S}$. 
The existence of the model structure follows directly from
Theorem \ref{Reedymodel1} except for the identification of
the Reedy cofibrations with the monomorphismsms.
For this, it suffices to  show that the acyclic v-fibrations  
are the trivial fibrations. But a map $f$ 
is an acyclic v-fibration iff
the map $\langle \delta_m\backslash f \rangle$
is trivial fibration for every $m\geq 0$ by \ref{Reedymodel1}
and by \ref{classicalmodel}. 
It then follows by \ref{trivfibbisimpset} that $f$ 
is an acyclic v-fibration iff it is a trivial fibration.
This completes the proof that the Reedy cofibrations are the monomorphismsms.
The model structure is left proper by \ref{Leftproper} since every 
object is cofibrant. 
Let us show that it is right proper. 
If $f:X\to Y$ is a v-fibration,
then the map $f_m:X_m\to Y_m$
is a Kan fibration for every $m\geq 0$ by \ref{Reedyfact2}.
Hence the base change
of a column-wise weak homotopy equivalence
along a v-fibration is a column-wise weak homotopy equivalence
since the model structure $({\cal C}_0, {\cal W}_0,{\cal F}_0)$
is proper.
Let us show that the model structure $({\cal C}'_0, {\cal W}'_0,{\cal F}'_0)$ 
is cartesian. See \ref{defcartclosedmod} for this notion.
By \ref{propQuillen2v}, it suffices to show that the functor $A\times (-):{\bf S}^{(2)}\to {\bf S}^{(2)}$ 
takes an acyclic map to an acyclic map. 
But this is clear since since the model structure $({\cal C}_0, {\cal W}_0,{\cal F}_0)$
is cartesian closed. Let us show that the model structure
is simplicial. 
The internal hom functor $(X,Y)\mapsto Y^X$
is a right Quillen functor of two variables (contravariant in the first) since the model
structure $({\cal C}'_0, {\cal W}'_0,{\cal F}'_0)$ 
is cartesian closed. 
Hence the composite $(X,Y)\mapsto Hom_2(X,Y)=i_2^*(Y^X)$
is also a right Quillen functor of two variables,
since the functor $i_2^*:{\bf S}^{(2)} \to {\bf S}$ 
is a right Quillen functor by \ref{Quillen pair for Reedy}.
\QED

\bigskip

Let $X$ be a bisimplicial set. From the map $[n]\to [0]$
we obtain a canonical map between the rows, $X_{\star 0}\to X _{\star n}$

\begin{definition}\label{defcategoricallyconstant} We shall say that a bisimplicial set $X$ {\it categorically constant}
if the canonical map
$X_{\star 0}\to X _{\star n}$ 
is a weak categorical equivalence for every $n\geq 0$.
\end{definition}

\begin{proposition}\label{categoricallyconstant} A v-fibrant bisimplicial set 
is categorically constant.  \end{proposition}

\noindent{\bf Proof}: 
If $i$ denotes the
inclusion $\Delta[0]\subseteq \Delta[n]$, then the map
$X/i:X/\Delta[n]\to X/\Delta[0]$ 
is a trivial fibration by \ref{Reedyfibrationvertically},
since $i$ is a monic weak homotopy equivalence.
But we have $(X/i)(X/t)=id$
since we have $ti=id$,  where $t$ is the map $\Delta[n]\to \Delta[0]$.
This shows by three-for-two that the map $X/t$
is a weak categorical equivalence.
\QED

We call a map of of bisimplicial sets $f:X\to Y$
a  {\it row-wise weak categorical equivalence} if the map $f_{\star n}:X_{\star n}\to Y_{\star n}$ 
is a weak categorical equivalence for every $n\geq 0$.

\begin{proposition}  \label{mapscategoricallyconstant} A map between v-fibrant simplicial sets $f:X\to Y$
is a row-wise weak categorical equivalence
iff it induces a weak categorical equivalence between the first rows. 
\end{proposition}

\noindent{\bf Proof}: 
If $f_{\star 0}: X_{\star 0}\to Y _{\star 0}$
is a weak categorical equivalence,
let us show that the map
$f_{\star n}: X_{\star n}\to Y _{\star n}$
is a weak categorical equivalence for every $n\geq 0$.
Consider the commutative square
$$
\xymatrix{
X_{\star 0} \ar[d]\ar[r]^{f_{\star 0}}  &Y_{\star 0}\ar[d] \\
X_{\star n} \ar[r]^{f_{\star n}} & Y_{\star n},}
$$
where the vertical maps are obtained from the map $\Delta[n]\to \Delta[0]$.
The vertical maps are  weak categorical equivalence
by \ref{categoricallyconstant}.
It follows by  three-for-two that 
the map $f_{\star n}$
is a weak categorical equivalence.
\QED

We shall say that $f:X\to Y$
is a  {\it horizontal fibration} or  
an {\it h-fibration} if the map
$\langle f/\delta_n \rangle $
is a quasi-fibration for every $n\geq 0$.
We shall say that a bisimplicial set $X$ is {\it h-fibrant}
if the map $X\to 1$ is an h-fibration.

\begin{proposition}  \label{ReedyQuasi-categoriesmodel} The category $ {\bf S}^{(2)}$ admits a 
model structure $({\cal C}^h_1, {\cal W}^h_1,{\cal F}^h_1)$ 
in which a cofibration is a monomorphism and 
a weak equivalence is a row-wise weak categorical equivalence.
A fibration is an h-fibration.
The model structure is left proper and cartesian closed.
The acyclic fibrations are the trivial  fibrations.
\end{proposition}

We call the model structure, the {\it horizontal model structure} on ${\bf S}^{(2)}$. 
For the notion of trivial fibration, see Definition \ref{trivfib}.

\medskip

\noindent{\bf Proof}: The model structure $({\cal C}^h_1, {\cal W}^h_1,{\cal F}^h_1)$  is the Reedy model structure 
$({\cal C}'_1, {\cal W}'_1,{\cal F}'_1)$
associated to the model structure 
$({\cal C}_1, {\cal W}_1,{\cal F}_1)$ of Theorem \ref{QCatmodel}.
It is similar to the vertical model structure
of Theorem \ref{verticalReedymodel}, except that
the weak equivalences are now defined row-wise
and by using weak categorical equivalences.
The existence of the model structure follows from Theorem \ref{Reedymodel1}.
Let us show that ${\cal C}'_1$ is the class of monomorphisms.
For this, it suffices to  show that the acyclic h-fibrations  
are the trivial fibrations. But a map $f$ 
is an acyclic h-fibration iff
the map $\langle f/\delta_n \rangle$
is trivial fibration for every $m\geq 0$ by \ref{Reedymodel1}
and by \ref{QCatmodel}. 
It then follows by \ref{trivfibbisimpset} that $f$ 
is an acyclic h-fibration iff it is a trivial fibration.
We have proved that ${\cal C}'_1$ is the class of monomorphisms.
It follows that every object is cofibrant. 
Hence the model structure is left proper by \ref{Leftproper}.
It remains to show that the model structure is cartesian closed.
By \ref{propQuillen2v},  it suffices to show that the functor $A\times (-):{\bf S}^{(2)}\to {\bf S}^{(2)}$ 
takes an acyclic map to an acyclic map. 
But this is clear since since the model structure $({\cal C}_1, {\cal W}_1,{\cal F}_1)$
is cartesian closed. \QED

\medskip

Recall the pair of adjoint functors $k_!: {\bf S} \leftrightarrow {\bf S}:k^!$
of Proposition \ref{Qpair3}.
By definition, we have $k_!(\Delta[n])= \Delta'[n]$ for every $n\geq 0$,
where $\Delta'[n]$ denotes
the (nerve of the) groupoid freely generated by
the category $[n]$.
Consider the functor $t:\Delta\times \Delta\to  {\bf S}$
defined by putting 
$t([m],[n])=\Delta[m]\times\Delta'[n]$
for every $m,n\geq 0$
and let $t_!: {\bf S}^{(2)}\to {\bf S}$
be the left Kan extension of the functor $t$ along the
Yoneda functor $\Delta^2\subset  {\bf S}^{(2)}$.
By definition, we have
$$t_!(\Delta[m]\Box \Delta[n])=\Delta[m]\times\Delta'[n].$$
The functor $t_!$ has a right adjoint $t^!: {\bf S} \to {\bf S}^{(2)}$.
If $X\in {\bf S}$, then 
$$t^!(X)_{mn}={\bf S}( \Delta[m]\times \Delta'[n],X)$$
for every $m,n\geq 0$.

\begin{lemma}  \label{Boxtoproductandt} There are 
natural
isomorphisms
$$t_!(A\Box B)=A\times k_!(B), \quad A\backslash t^!(X)=k^!(X^A)
 \quad {\rm and}\quad  t^!(X)/B=X^{k_!(B)}$$
for $A,B\in {\bf S}$ and $X\in {\bf S}^{(2)}$.
\end{lemma}

\noindent{\bf Proof}: The functor $(A,B)\mapsto t_!(A\Box B)$
is cocontinuous in each variable and it extends
the functor 
$([m],[n])\mapsto \Delta[m]\times \Delta'[n]$
along the 
Yoneda functor $\Delta\times \Delta\to  {\bf S}^{(2)}.$
Similarly for the functor $(A,B)\mapsto A\times k_!(B)$,
since $k_!(\Delta[n])= \Delta'[n]$.
It follows that there is a natural isomorphism
$t_!(A\Box B)=A\times k_!(B)$.
The first statement of the proposition is proved.
Let us prove the second.
The functor $X\mapsto A\backslash t^!(X)$
is right adjoint to the functor
$B\mapsto t_!(A\Box B)$, since a composite
of right adjoints is right adjoint to the composite in reverse order.
Similarly, the functor $X\mapsto k^!(X^A)$
is right adjoint to the functor $B\mapsto A\times k_!(B)$. 
This proves the result by uniqueness of adjoint 
since $t_!(A\Box B)=A\times k_!(B)$.
Let us prove the third. 
The functor $X\mapsto  t^!(X)/B$
is right adjoint to the functor
$A\mapsto t_!(A\Box B)$, since a composite
of right adjoints is right adjoint to the composite in reverse order.
Similarly, the functor $X\mapsto X^{k_!(B)}$
is right adjoint to the functor $A\mapsto A\times k_!(B)$. 
This proves the result by uniqueness of adjoints 
since $t_!(A\Box B)=A\times k_!(B)$.
\QED

\begin{theorem}  \label{Qpairforhorizontalandvertical} 
The adjoint pair of functors 
$$t_!:{\bf S}^{(2)} \leftrightarrow {\bf S}:t^!$$
is a Quillen pair
between the horizontal (resp. vertical) model structure on ${\bf S}^{(2)}$
and the model structure for quasi-categories on $ {\bf S}$.
\end{theorem}

\noindent{\bf Proof}: Let us first show that $(t_!,t^!)$ is
a Quillen pair between the horizontal
model structure
and the model structure for quasi-categories.
We shall use
the criteria of \ref{Quillenpair2}.
We first verify that $t_!$
takes a cofibration to a cofibration.
Let ${\cal A}$ be the class of maps
$u\in  {\bf S}^{(2)} $ such that
$t_!(u)$ is monic. The class  ${\cal A}$ is saturated
since the functor $t$ is cocontinuous.
Let us show that the map $\delta_m\Box' \delta_n$
belongs to ${\cal A}$ for every $m,n\geq 0$.
We have $t_!(\delta_m\Box' \delta_n)=\delta_m\times' k_!(\delta_n).$
by Lemma \ref{Boxtoproductandt}.
But $k_!(\delta_n)$ is monic since the functor $k_!$
takes a monomorphism to a monomorphism by
\ref{Qpair3}. Hence the map $\delta_m\times' k_!(\delta_n)$
is monic and this shows that the map $\delta_m\Box' \delta_n$
belongs to ${\cal A}$. It follows by \ref{saturatedmonobisset}
that every monomorphism
belongs to ${\cal A}$.
We have proved that $t_!$
takes a cofibration to a cofibration.
Let us now show that $t^!$
takes a fibration 
to a fibration.
For this we have to show that
if $f:X\to Y$ is a quasi-fibration 
then the map $t^!(f):t^!(X)\to t^!(Y)$
is an h-fibration. For this
it suffices to show that the map 
$\langle t^!(f)/u\rangle $ is a quasi-fibration for every monomorphism of simplicial
sets $u:A\to B$. 
But the square
$$\xymatrix{
 t^!(X)/B \ar[r]\ar[d]&  t^!(X)/A \ar[d] \\
 t^!(Y)/B   \ar[r]&  t^!(Y)/A
 }
 $$
is isomorphic to the square 
$$\xymatrix{
 X^{k_!(B)} \ar[r]\ar[d]& X^{k_!(A)} \ar[d] \\
 Y^{k_!(B)}   \ar[r]&  Y^{k_!(A)}
 }
 $$
by Lemma \ref{Boxtoproductandt}.
Hence the map $\langle  t^!(f)/u\rangle$
is isomorphic to the map
$\langle k_!(u),f\rangle$.
But the map $k_!(u)$ is monic. 
Hence the $\langle k_!(u),f\rangle$
is a quasi-fibration
by \ref{QCatmodel} since $f$
is a quasi-fibration.
We have proved $(t_!,t^!)$
is a Quillen pair between the horizontal
model structure
and the model structure for quasi-categories.
Let us now show that it
is a Quillen pair between the vertical
model structure
and the model structure for quasi-categories.
We saw that $t_!$
takes a cofibration to a cofibration.
Let us now show that $t^!$
takes a fibration to a fibration.
For this we have to show that
if $f:X\to Y$ is a quasi-fibration
then the map $t^!(f):t^!(X)\to t^!(Y)$
is a v-fibration. 
For this it suffices to show that the map 
$\langle u\backslash t^!(X)\rangle $ is a Kan fibration for every monomorphism of simplicial
sets $u:A\to B$. 
But the square
$$\xymatrix{
B\backslash t^!(X) \ar[r]\ar[d]& A\backslash t^!(X) \ar[d] \\
B\backslash t^!(Y)   \ar[r]& A\backslash t^!(Y)
 }
 $$
is isomorphic to the square 
$$\xymatrix{
 k^!(X^B) \ar[r]\ar[d]& k^!(X^A ) \ar[d] \\
 k^!(Y^B)   \ar[r]&  k^!(Y^A)
 }
 $$
by Lemma \ref{Boxtoproductandt}.
Hence the map $\langle u\backslash t^!(X)\rangle$
is isomorphic to the map
$k^!\langle u,f\rangle$.
 But $\langle u,f\rangle$ is a quasi-fibration
by Theorem \ref{QCatmodel}. 
Thus, $k^!\langle u,f\rangle$
is a Kan fibration by \ref{Qpair3}.
We have proved that $(t_!,t^!)$
is a Quillen pair between the vertical
model structure
and the model structure for quasi-categories.
\QED

\section {Segal spaces}

Segal spaces were introduced by Rezk in \cite{Rez}.
They are the fibrant objects of a model structure
on the category of bisimplicial sets.
The model structure is a Bousfield
localisation of the (vertical) Reedy model structure introduced in the previous section.

\bigskip

Let  $I_n\subseteq \Delta[n]$ be the $n$-chain. 
For any simplicial space $X$ we have a canonical bijection
$$
I_n \backslash X =X_1  \times_{\partial_0,\partial_1} X_1  \times \cdots 
 \times_{\partial_0,\partial_1}  X_1, 
$$
where the successive  fiber products are calculated by using the face maps
$\partial_0,\partial_1:X_1\to X_0$.

\begin{definition}  \cite {Rez} \label{defSegalcondition} We shall say that a simplicial space $X$ 
satisfies the  {\it Segal condition} 
if the map 
$$i_n\backslash  X :\Delta[n]\backslash  X \longrightarrow I_n \backslash X$$
obtained from the inclusion $i_n:I_n\subseteq \Delta[n]$
is a weak homotopy equivalence for every $n\geq 2$.
A {\it Segal space}  is a $v$-fibrant simplicial space which satisfies the Segal condition,
\end{definition}

\medskip

We shall say that a map of simplicial
spaces $u:A\to B$ is a {\it Segal weak equivalence} if
the map
$$Hom_2(u,X):Hom_2(B,X)\to Hom_2(A,X)$$
is a weak homotopy equivalence
for every Segal space $X$.

\bigskip

\begin{theorem}  \label{Segalspacemodel} \cite{Rez} The category
$({\bf S}^{(2)},Hom_2)$ admits a simplicial model structure
$({\cal C}_S, {\cal W}_S,{\cal F}_S)$ in which  ${\cal C}_S$ is the class of monomorphisms and 
${\cal W}_S$ is the class of Segal weak equivalences.
The model structure is left proper and cartesian closed.
The acyclic fibrations are the trivial  fibrations.
The fibrant objects are the Segal spaces.
\end{theorem}

We shall say that $({\cal C}_S, {\cal W}_S,{\cal F}_S)$ is  the {\it model 
structure for Segal spaces}. 
The model structure 
is a Bousfield localisation of the vertical model structure 
$({\cal C}^v_0, {\cal W}^v_0,{\cal F}^v_0)$ 
of Theorem \ref{verticalReedymodel}.
For the notion of Bousfield localisation, see Definition \ref{DefBousfieldLocalisationA}. 
In particular, a map between Segal spaces
is acyclic (resp. a fibration)  iff it is a column-wise weak homotopy equivalence
(resp. a v-fibration).

\medskip

Recall the total space functor $t_!:{\bf S}^{(2)} \to {\bf S}$
of Theorem \ref{Qpairforhorizontalandvertical}.

\begin{theorem}  \label{QpairforSegalspaces} 
The adjoint pair of functors 
$$t_!:{\bf S}^{(2)} \leftrightarrow {\bf S}:t^!$$
is a Quillen pair
between the model category for Segal spaces 
and the model category for quasi-categories.
\end{theorem}

\noindent{\bf Proof}: We saw in \ref{Qpairforhorizontalandvertical} 
that the pair of adjoint functors $(t_!,t^!)$ is a 
Quillen pair
between the vertical model category 
and the model category for quasi-categories.
Hence it suffices to show by \ref{Quillenpair2}.
that
the functor $t^!$ takes a quasi-category to a Segal space.
If $X$ is a quasi-category then $t^!(X)$ 
is v-fibrant by \ref{Qpairforhorizontalandvertical}.
Let us show that $t^!(X)$ 
satisfies the Segal condition.
For this, it suffices to show that
the map $i_n\backslash t^!(X)$
is a trivial fibration for every $n\geq 0$, where $i_n$
denotes the inclusion $I_n\subseteq \Delta[n]$.
But the map $i_n\backslash t^!(X)$
is isomorphic to the map $k^!(X^{i_n})$
by Lemma \ref{Boxtoproductandt}. 
The map $i_n$ is weak categorical equivalence by \ref{midcofibration}
since it  is mid anodyne by \ref{midado2}.
Hence the map
$X^{i_n}$ is a trivial fibration 
by \ref{QCatmodel} since $i_n$ is monic.
It follows that the map $k^!(X^{i_n})$ 
is a trivial fibration since $k^*$ is a right Quillen functor
by \ref{Qpair3}. We have proved that  the map
$i_n\backslash t^!(X)$
is a trivial fibration.
This shows that  $t^!(X)$ is a Segal space.
\QED

\bigskip

\begin{proposition}\label{Segalspace2} Let $X$ be a v-fibrant simplicial space. Then
the following conditions are equivalent:
\begin{itemize}
\item {\rm (i)} $X$ is a Segal space 
\item {\rm (ii)} The map $h^k_n\backslash X$ is a trivial fibration
for every $0<k<n$;
\item {\rm (iii)} the map
$u \backslash X$
is a trivial fibration
for every mid anodyne map $u\in  {\bf S}$;
\item {\rm (iv)} the map
$X\slash \delta_n$
is a mid fibration for every $n\geq 0$;
\item {\rm (v)}
the map
$X\slash v$
is a mid fibration for every monomorphism $v\in {\bf S}$
\end{itemize}
\end{proposition}

The notion of mid anodyne map was defined in \ref{mid cofibration def}.  The proof depends on the following lemma. 

\medskip

We shall say that a class of  monomorphisms ${\cal  A}\subseteq {\bf S}$
has the {\it right cancellation property} if the implication 
$$vu\in  {\cal A }\quad {\rm and} \quad  u \in  {\cal A }\quad \Rightarrow \quad v\in  {\cal A }$$
is true for any pair of monomorphisms $u: A  \to B $
and $v:B\to C$.

\medskip

Let  $I_n\subseteq \Delta[n]$ be the $n$-chain.

\begin{lemma}\label{midado3}
Let ${\cal A}\subseteq {\bf S}$ be a saturated class of monomorphisms
having the right cancellation property.
If the inclusion $I_n\subseteq \Delta[n]$ belongs to ${\cal A }$ for every $n\geq 2$,
then every  mid anodyne map belongs to ${\cal A }$.
\end{lemma}

\noindent {\bf Proof}: 
We shall first prove that 
the inclusion $I_n\subseteq \partial_0\Delta[n] \cup \partial_n\Delta[n]$
belongs to ${\cal A }$ for every $n>1$.
This is obvious if $n=2$ since we have $I_n= \partial_0\Delta[n] \cup \partial_n\Delta[n]$
 in this case. Let us suppose $n>2$.
It suffices to show  
that each inclusion
  $$\xymatrix{ I_n\  \ar@{^{(}->}[r]^(0.35){i_n} &\  I_n    \cup   \partial_n\Delta[n] \ \ar@{^{(}->}[r]^(0.45){j_n}& \partial_0\Delta[n] \cup \partial_n\Delta[n] 
  }$$
 belongs to ${\cal A }$, since ${\cal A }$ is closed under composition.
 The square of inclusionsÊ
 $$
\xymatrix{
I_{n-1} \ar[d]\ar[r]&  I_n   \ar[d]^{i_n}\\
\partial_n\Delta[n]  \ar[r]  &\partial_n\Delta[n] \cup    I_n }.
$$
is a pushout since $\partial_n\Delta[n] \cap   I_n =I_{n-1}$.
Thus, $i_n\in {\cal A }$ since ${\cal A }$
is closed under cobase change and
since the inclusion $I_{n-1}\subset \Delta[n-1]=\partial_n\Delta[n]$
belongs to $ {\cal A }$.
It remains to show that the inclusion $j_n$ belongs to ${\cal A }$.
If $d_0:\Delta[n-1] \to \Delta[n]$, then we have
$d_0^{-1}(I_n)=I_{n -1}$ 
and  $d_0^{-1}(\partial_n\Delta[n])=\partial_{n-1}\Delta[n-1] $.
It follows from this observation that the following square is a
pushout,
$$
\xymatrix{
 I_{n -1} \cup \partial_{n-1}\Delta[n-1] \ar[d]^{k_{n-1}}\ar[r] & I_n    \cup   \partial_n\Delta[n]\ar[d]^{j_n}\\
\Delta[n-1]  \ar[r]^(0.40){d'_0}	 &\partial_0\Delta[n] \cup \partial_n\Delta[n],}
$$
where $k_{n-1}$ is the inclusion and where $d'_0$ is induced by $d_0$.
Let us show that $k_{n-1}\in {\cal A }$.
The composite 
$$\xymatrix{ I_{n-1}\  \ar@{^{(}->}[rr]^(0.35){i_{n-1}} && I_{n -1} \cup \partial_{n-1}\Delta[n-1]  \ \ar@{^{(}->}[rr]^(0.65){k_{n-1}}&& \Delta[n-1]  
  }$$
belongs to ${\cal A }$ by assumption.
We have $i_{n-1}\in {\cal A }$
by the same argument as above
since the inclusion $I_{n-2}\subset  \Delta[n-2]$
belongs to $ {\cal A }$.
It follows by the right cancellation property of the class ${\cal A }$
that $k_{n-1}$ belongs to ${\cal A }$.
Thus, $j_n\in {\cal A }$ since the class ${\cal A }$
is closed under cobase change.
This completes the proof that
the inclusion $I_n\subset \partial_0\Delta[n] \cup \partial_n\Delta[n]$
 belongs to ${\cal A }$ for $n>2$, hence also for $n>1$.
We can now prove the lemma. 
For this
it suffices to show that
the inclusion  $ \Lambda^k[n]  \subset \Delta[n]$
belongs to the class
${\cal A }$ for every $0<k<n $, since the class ${\cal A }$ is saturated.
This is obvious if $n=2$
since $\Lambda^1[2] =I_2$.
We shall suppose $n>2$.
By the cancellation property of the class $\cal A $, it suffices to show that
the inclusion  
$I_n\subseteq \Lambda^k[n] $
belongs to ${\cal A }$,
since the inclusion
$I_n\subseteq \Delta[n] $
belongs to ${\cal A }$.
If $S$ is a subset of $[n]$, let us put
$$\Lambda^S[n] =\bigcup_{i\not\in S} \partial_i\Delta[n].$$
We shall prove
that if $n>2$ and $S$
is a non-empty subset of the interval $[1,n-1]$, then
the inclusion  
$I_n\subset  \Lambda^S[n]$
belongs to ${\cal A }$.
We argue by induction
on $n>2$  
and $s=n-{\rm Card}(S)>0$.
If $s=1$, then $S=[1,n-1]$ and $\Lambda^S[n]= \partial_0\Delta[n] \cup \partial_n\Delta[n] $.
The result was proved above in this case. 
If $s >1$ let us choose an element $b\in [1,n-1]\setminus S$
and put $T=S\cup \{b\}$.
The inclusion  
$I_n\subset  \Lambda^T[n] $
belongs to ${\cal A }$
by the induction hypothesis
 since $n-{\rm Card}(T )<s$.
Let us show that the  inclusion $ \Lambda^T [n]\subset \Lambda^S [n]$
belongs to ${\cal A }$.
The square 
$$
\xymatrix{
\partial_{b}\Delta[n]\cap\Lambda^T [n] \ar[d] \ar[r]&  \Lambda^T [n]  \ar[d]\\
\partial_{b}\Delta[n] \ar[r]   &\Lambda^S [n]   }
$$
is a pushout since $\Lambda^S [n] =\partial_b\Delta[n]\cup \Lambda^T[n]$.
Hence it suffices to show 
that the inclusion $\partial_{b}\Delta[n]\cap\Lambda^T[n]\subset \partial_{b}\Delta[n]$
belongs to ${\cal A }$.
Let
$U\subseteq [n-1]$
be the inverse image of the subset $T$ by the map $d_b:[n-1]\to [n]$.
The inclusion $\partial_{b }\Delta[n]\cap\Lambda^T [n]\subset \partial_{b}\Delta[n]$
is isomorphic to the inclusion
$\Lambda^U [n-1]\subset \Delta[n-1]$.
Hence it suffices to show that the latter 
belongs to ${\cal A }$.
The subset $U$ is non-empty since it is in bijection with $S$.
Moreover, $U\subseteq [1,n-2]$
since $S\subseteq [1,n-1]$ and $0<b<n$.
Hence
the inclusion $I_{n-1}\subset \Lambda^U [n-1]$
belongs to ${\cal A }$ by the induction hypothesis on $n$.
It follows that the inclusion
$\Lambda^U [n-1]\subset \Delta[n-1]$ belongs to ${\cal A }$
by the cancellation property of the class ${\cal A }$. Hence
the inclusion $ \Lambda^T [n]\subset \Lambda^S [n]$
belongs to ${\cal A }$, since  ${\cal A }$ is closed under cobase change.
It then follows by composing that the inclusion $I_n\subset  \Lambda^S[n] $
belongs to ${\cal A }$, since the  inclusion  
$I_n\subset  \Lambda^T[n] $
belongs to ${\cal A }$
by the induction hypothesis.
We have proved
that if $S$
is a non-empty subset of the interval $[1,n-1]$, then
the inclusion  
$I_n\subset  \Lambda^S[n]$
belongs to ${\cal A }$.
In particular, this shows that 
the inclusion  
$I_n\subseteq \Lambda^k[n] $
belongs to ${\cal A }$.
\QED

\noindent
{\bf Proof of proposition \ref{Segalspace2}}: 
Let us prove the implication
(i)$\Rightarrow$(iii).
Let 
${\cal A }$ be the class of monomorphisms $u\in {\bf S}$
such that the map $u \backslash X$
is a weak homotopy equivalence.
It is obvious from this definition
that ${\cal A }$ has the right cancellation property.
Let us verify that ${\cal A }$
is saturated.
A monomorphism $u\in {\bf S}$ belongs to ${\cal A }$
iff the map $u \backslash X$ is a trivial fibration
since it is a Kan fibration by \ref{Reedyfibrationvertically}.
But $u \backslash X$ is a trivial fibration 
iff we have $\delta_n\pitchfork (u\backslash X)$  for every $n\geq 0$
by \ref{saturatedmonosset}.
But the  condition 
$\delta_n\pitchfork (u\backslash X)$  
is equivalent to the condition 
$u\pitchfork  (X\slash \delta_n)$
by  \ref{lift2}.
Thus, a monomorphism $u$ belongs to ${\cal A}$
iff we have $u\pitchfork (X\slash \delta_n)$
for every $n\geq 0$.
This shows that the class ${\cal A }$
is saturated. If $X$ is a Segal space then 
the inclusion $I_n\subseteq \Delta[n]$
belongs to ${\cal A}$ for every $n\geq 2$.
It then follows from \ref{midado3} that
every mid anodyne map belongs to
${\cal A}$. The implication (i)$\Rightarrow$(iii)
is proved. 
The converse
follows from the fact that the inclusion $I_n\subseteq \Delta[n]$
is mid anodyne by \ref{midado2}.
The implication (iii)$\Rightarrow$(ii)
is obvious. Let us prove the implication (ii)$\Rightarrow$(v).
If $v\in  {\bf S}$ is monic, let us show that
$X\slash v $
is a mid fibration.
But the condition
$h^k_n \pitchfork (X/v)$
is equivalent to the condition $v\pitchfork (h^k_n\backslash  X)$
by \ref{lift2}. This shows that we have $h^k_n \pitchfork (X/v)$
for every $0<k<n$, since the map $h^k_n\backslash  X$
is a trivial fibration in this case.
The implication (ii)$\Rightarrow$(v)
is proved. 
The implication (v)$\Rightarrow$(iv)
is obvious. 
The implication (iv)$\Rightarrow$(iii)
follows from \ref{Reedyfib1}.
\QED

\begin{corollary} \label{Segalspacetoquasi} If $X$ is a Segal space, then the simplicial set $X/A$ is a 
quasi-category for any simplicial set $A$. In particular, every row of $X$ is a quasi-category.
\end{corollary}

\noindent {\bf Proof}: 
 If $i_A$ denotes the inclusion $\emptyset \subset  A$,
then the map $X/i_A: X/ A\to X/ \emptyset$
is a mid fibration
by (v). This shows that $X/A$
is a quasi-category since
$X/ \emptyset=1$.
In particular, the $n$th row $X_{\star n}=X/\Delta[n]$
is a quasi-category.
\QED

\begin{lemma}\label{Anodyne}
Let ${\cal A}\subseteq {\bf S}$ be a saturated class of monomorphisms
having the right cancellation property.
If every face map $d_i:\Delta[n-1]\subset \Delta[n]$ belongs to ${\cal A }$,
then every anodyne map belongs to ${\cal A }$.
\end{lemma}

See \ref{anodyne def} for the notion of anodyne map.
\medskip

\noindent{\bf Proof}: It suffices to show that
every horn $h^k_n:\Lambda^k[n]\subset \Delta[n]$
belongs to ${\cal A}$, since ${\cal A}$ is saturated.
More generally, if $S$ is a proper non-empty subset of $[n]$, let us put
$$\Lambda^S[n] =\bigcup_{i\not\in S} \partial_i\Delta[n].$$
We shall prove by induction on $n\geq 1$
that the inclusion $\Lambda^S[n] \subset \Delta[n]$
belongs to ${\cal A}$.
The result is clear if $n=1$, since  $h^0_1=d_1$
and $h^1_1=d_0$. Let us suppose $n>1$.
It suffices to show that 
the inclusion $\partial_i\Delta[n]\subset \Lambda^S[n]$
belongs to ${\cal A}$ for some $i\not\in  S$, since the
class ${\cal A}$ has the  right cancellation property.
We have $\Lambda^S[n]=\partial_i\Delta[n]$ if $S=[n]\setminus\{i\}$.
Hence it suffices to show
that the inclusion $\Lambda^T[n]\subseteq \Lambda^S[n]$
belongs to ${\cal A}$ for any pair of proper non-empty subsets
$S\subset T\subset [n]$. 
For this it suffices to consider the case where $T=S\cup \{t\}$ with $t\not\in S$,
since the class ${\cal A}$ is closed under composition.
The square 
$$
\xymatrix{
\partial_{t}\Delta[n]\cap\Lambda^T [n] \ar[d] \ar[r]&  \Lambda^T [n]  \ar[d]\\
\partial_{t}\Delta[n] \ar[r]   &\Lambda^S [n]   }
$$
is a pushout since $\Lambda^S [n] =\partial_t\Delta[n]\cup \Lambda^T[n]$.
Hence it suffices to show 
that the inclusion $\partial_{t}\Delta[n]\cap\Lambda^T[n]\subset \partial_{t}\Delta[n]$
belongs to ${\cal A }$.
Let
$U\subseteq [n-1]$
be the inverse image of the subset $T$ by the map $d_t:[n-1]\to [n]$.
The inclusion $\partial_{t}\Delta[n]\cap\Lambda^T [n]\subset \partial_{t}\Delta[n]$
is isomorphic to the inclusion
$\Lambda^U [n-1]\subset \Delta[n-1]$.
Hence it suffices to show that the latter 
belongs to ${\cal A }$.
The subset $U$ is non-empty since it is in bijection with $S$.
Moreover $U$ is proper, since ${\rm Card}(U)={\rm Card}(S)<{\rm Card}(T)$
and $T$ is a proper subset of $[n]$.
Thus, the inclusion
$\Lambda^U [n-1]\subset \Delta[n-1]$
belongs to ${\cal A }$ by the induction
hypothesis. 
This proves that the inclusion $\Lambda^T[n]\subseteq \Lambda^S[n]$
belongs to ${\cal A }$. \QED

\begin{lemma} \label{Segalspacelemma3} A mid fibration between quasi-categories is a trivial fibration iff it is a weak categorical equivalence surjective on vertices.
\end{lemma}

\noindent{\bf Proof}: The necessity is clear. Conversely, if a mid fibration between quasi-categories $f:X\to Y$ is a weak categorical equivalence surjective on vertices,
let us show that it is a trivial fibration. For this it suffices to show that $f$ is a quasi-fibration.
But for this, it suffices to show that the functor $\tau_1(f)$ is a quasi-fibration by \ref{qfbetqc}, since
$f$ is a map between quasi-categories.
But the functor $\tau_1(f)$ is an equivalence of categories, since $f$ is a weak categorical equivalence by assumption.
Moreover, the functor $\tau_1(f)$ is a surjective on objects,
since the map $f$ is surjective on vertices by assumption.
Thus, $\tau_1(f)$ is an acyclic fibration for the natural model structure on $\bf Cat$.
It follows that $\tau_1(f)$ is a quasi-fibration, and hence that $f$
is a trivial fibration.
\QED

\begin{proposition} \label{Segalspace3} A bisimplicial set $X$
is a Segal space
iff the following three conditions are satisfied:
\begin{itemize}
\item {\rm (i)} the map
$X\slash \delta_n $
is a mid fibration for every $n\geq 0$;
\item {\rm (ii)} $X_0$ is a Kan complex;
\item {\rm (iii)} $X$ is categorically constant.
\end{itemize}
\end{proposition}

\noindent{\bf Proof}: ($\Rightarrow$)
The simplicial set $X_0$ is a Kan  complex by \ref{Reedyfact2}, since $X$ is $v$-fibrant by
assumption.
Moreover, $X$ is categorically constant
by \ref{categoricallyconstant}.Ê
The map $X\slash \delta_n$
is a mid fibration for every $n\geq 0$ by Proposition \ref{Segalspace2}
since $X$ is a Segal space.
($\Leftarrow$) We shall prove that $X$ is vertically fibrant.
by showing  
that the map $X/v$ is a trivial fibration for every
anodyne map $v\in {\bf S}$.
Observe first that the map $X\slash v$
is a mid fibration for every monomorphism $v:S\to T$ by \ref{Reedyfib1},
since the map $X\slash \delta_n $
is a mid fibration for every $n\geq 0$ by assumption.
In particular, 
$X/S$ is a quasi-category
for any simplicial set $S$, since the map
$X/S\to X/ \emptyset=1$
is a mid fibration.
Hence the map $X\slash v $
is a mid fibration between quasi-categories for any monomorphism $v:S\to T$.
We claim that if $v$ is anodyne and $X\slash v$ is a weak categorical equivalence,
then $X\slash v$ is actually a trivial fibration.
For this, it suffices to show
that $X\slash v$ is surjective on vertices by Lemma \ref{Segalspacelemma3}.
But every map $S\to X_0$ can be extended along $v$ to a map $T\to X_0$,
since $v$ is anodyne and $X_0$ is a Kan complex by assumption.
Hence the map $(X/v)_0:(X/T)_0\to (X/S)_0$
is surjective on vertices. 
This shows that $X\slash v$ is a trivial fibration if $v$ is anodyne and $X\slash v$ is a weak categorical equivalence.
Let us now prove that $X$ is vertically fibrant.
For this  it suffices to show 
that the map $X/v$ is a trivial fibration for every
anodyne map $v\in {\bf S}$ by \ref{Reedyfibrationvertically}. 
Let ${\cal A}$ be the class of anodyne maps
$v\in {\bf S}$ such that the map $X/v$ is a trivial fibration.
We shall prove that every anodyne map belongs to ${\cal A}$ by using Lemma \ref{Anodyne}.
An anodyne map belongs 
${\cal A}$ iff the map $X/v$ is a weak categorical equivalence by the above.
It follows by three-for-two that ${\cal A}$ has the right cancellation property.
Let us show that ${\cal A}$ is saturated.
The condition $\delta_n\pitchfork (X\slash v)$
is equivalent to the condition $v\pitchfork (\delta_n\backslash X)$
by \ref{lift2}. Thus, an anodyne map
$v$ belongs to ${\cal A}$ iff we have $v\pitchfork (\delta_n\backslash X)$
for every $n\geq 0$. This description implies that
the class ${\cal A}$ is saturated.
Let us show that every face map $d_i:\Delta[n-1]\to \Delta[n]$ belongs to  ${\cal A}$.
The canonical maps $X_{\star 0}\to X_{\star n}$ and $X_{\star 0}\to X_{\star ,n-1}$ are weak categorical equivalences, since $X$
is categorically constant. It follows by three-for-two that the map $X/d_i:X_{\star n}\to X_{\star ,n-1}$
is a weak categorical equivalence. This proves that $d_i\in {\cal A}$, since $d_i$
is anodyne.
It then follows by Lemma \ref{Anodyne} 
that every anodyne map belongs to $\cal A$.
Thus, $X$ is v-fibrant by \ref{Reedyfibrationvertically}.
Let us now show that $X$ satisfies the Segal condition.
For this, it suffices to show that the map $ i_m\backslash X$
is a trivial fibration for every $m\geq 0$,
where $i_m$ denotes the inclusion 
of the $m$-chain $I_m\subseteq \Delta[m]$. But 
for this it suffices to show that we
have $\delta_n \pitchfork (i_m\backslash X)$ for every $n\geq 0$.
But $i_m$
is mid anodyne by Lemma \ref{midado2}. 
Hence we have
$i_m\pitchfork (X\slash \delta_n) $ by (i). 
It follows that we have $\delta_n \pitchfork (i_m\backslash X)$ 
by \ref{lift2}. 
This shows that $X$ is a Segal space.
\QED

\begin{proposition}\label{v-fibSegalspace}
Let $f: X \ra Y$ be a v-fibration between Segal spaces.  
Then the map
$$\langle u\backslash f\rangle :B \backslash X \longrightarrow B \backslash Y \times_{A \backslash Y  }A\backslash X $$
is a trivial fibration 
for any mid anodyne map $u:A\to B$. 
Moreover, the map
$$\langle f/v \rangle :X/T\longrightarrow Y/T\times_{Y/S}X/S$$
is a mid fibration between quasi-categories
for any monomorphism of simplicial sets
$v:S\to T$.
\end{proposition}

\noindent
{\bf Proof}: Let us prove the first statement. The map $\langle u\backslash f\rangle$
 is a Kan fibration by \ref{Reedyfibrationvertically}. 
 Let us show
 that it is a weak homotopy equivalence.
The horizontal maps in the commutative square 
 $$
\xymatrix{
B \backslash X  \ar[d] \ar[r]& A\backslash X  \ar[d]\\
  B \backslash Y \ar[r]  & A \backslash Y     }
$$
are
trivial fibrations by \ref{Segalspace2} .
It follows that 
$\langle u\backslash f\rangle$
is a weak homotopy equivalence.
It is thus a trivial fibration since it is a Kan fibration.
Let us prove the second statement. 
The domain of $\langle f/v \rangle$ 
is a quasi-category by \ref{Segalspace2}.
Let us show that its codomain is a quasi-category.
The projection $p_2$ in the pullback square
$$
\xymatrix{
Y/S \times_{Y/S}X/T  \ar[d] \ar[r]^(0.65){p_2}& X/S  \ar[d]\\
 Y/T \ar[r]  & Y/S   }
$$
is a mid fibration since the
bottom map is a mid fibration by \ref{Segalspace2}.
It follows that the domain of $p_2$ 
is a quasi-category since its codomain
is a quasi-category by \ref{Segalspace2}.
Let us now show that
$\langle f/v \rangle$ is a mid fibration.
By \ref{midcofibrationfibration} it suffices to show that we have
$u\pitchfork \langle f/v\rangle$
for every mid anodyne map $u:A\to B $.
But the condition  $u\pitchfork \langle f/v\rangle$
is equivalent to the condition
$v\pitchfork \langle u\backslash f\rangle$
by \ref{lift2}.
But we have $v\pitchfork \langle u\backslash f\rangle$
since 
$\langle u\backslash f\rangle$ is a trivial fibration by
the first part of the proof.
\QED

\section {Two equivalences with complete Segal spaces}

Complete Segal spaces were introduced by Charles Rezk in \cite{Rez}.
They are the fibrant objects of 
 a Quillen model structure on the category of 
 simplicial spaces $[\Delta^o, {\bf S}]= {\bf S}^{(2)}$ . 
 We call this model structure the {\it model 
structure for complete Segal spaces}. 
The goal of the section is to
establish two Quillen equivalences 
$$p_1^*:{\bf S}\leftrightarrow {\bf S}^{(2)}:i_1^* \quad {\rm and} \quad t_!:{\bf S}^{(2)} \leftrightarrow {\bf S}:t^!.$$
 between the model category for quasi-categories
and the model category for complete Segal spaces.
The functor $i_1^*$  associates to a bisimplicial set $X$
its first row  $X_{\star 0}$. It shows that a complete Segal 
space is determined up to equivalence by its first row.
The functor $t_!$  associates to a bisimplicial set $X$
a total simplicial set $t_!X$.

\bigskip

Let $J$ be the groupoid generated by
one isomorphism $0\to 1$.
We shall regard $J$ as a simplicial set via the nerve functor.
A Segal space $X$
is said to be {\em complete } if 
the map
$$1\backslash X  \longrightarrow J \backslash X$$
obtained from the map $J\to 1$
is a weak homotopy equivalence.

\medskip

We shall say that a map of simplicial
spaces $u:A\to B$ is a {\it Rezk weak equivalence} if
the map
$$Hom_2(u,X):Hom_2(B,X)\to Hom_2(A,X)$$
is a weak homotopy equivalence
for every complete Segal space $X$.

\bigskip

\begin{theorem} [Rezk, see  \cite{Rez}]  \label{Rezkmodel} The simplicial category $({\bf S}^{(2)},Hom_2)$ admits a
simplicial model structure
$({\cal C}_R, {\cal W}_R,{\cal F}_R)$ in which  ${\cal C}_R$ is the class of monomorphisms and 
${\cal W}_R$ is the class of Rezk weak equivalences.
The model structure is left proper and cartesian closed.
The acyclic fibrations are the trivial  fibrations.
The fibrant objects are the complete Segal spaces.
\end{theorem}

We call $({\cal C}_R, {\cal W}_R,{\cal F}_R)$ the {\it Rezk model 
structure} or the {\it model 
structure for complete Segal spaces}. 
For the notion of trivial fibration, see Definition \ref{trivfib}.
The model structure 
is a Bousfield localisation of the Segal space model structure 
of Theorem \ref{Segalspacemodel}. Hence it is also 
a Bousfield localisation of the vertical model structure 
$({\cal C}^v_0, {\cal W}^v_0,{\cal F}^v_0)$ 
of Theorem \ref{verticalReedymodel}.
For the notion of Bousfield localisation, see Definition \ref{DefBousfieldLocalisationA}. 
In particular, a map between complete Segal spaces
is acyclic (resp. a fibration)  iff it is a column-wise weak homotopy equivalence
(resp. a v-fibration).

\medskip

The first projection $p_1: \Delta\times \Delta \to \Delta$
is left adjoint to the functor $i_1:\Delta \to  \Delta\times \Delta$ 
defined by putting 
$i_1([n])=([n],[0])$ for every $n\geq 0$. 
We thus obtain a pair of adjoint
functors 
$$p_1^*:{\bf S}\leftrightarrow {\bf S}^{(2)}:i_1^*.$$
If $X$ is a bisimplicial set, then $i_1^*(X)$
is the first row of $X$. 
Notice that we have $p_1^*(A)=A\Box 1$ for every simplicial set $A$.
In \ref{Theorem A} we shall prove that

\bigskip

\noindent{\bf Theorem} \emph{
The adjoint pair of functors $$p^*_1: {\bf S}\leftrightarrow  {\bf S}^{(2)}:i_1^*$$ is a Quillen equivalence
between the model category for  quasi-categories
and the model category for complete Segal spaces}.

\bigskip

Recall the ``total space" functor 
$t_!:{\bf S}^{(2)} \to {\bf S}$
of Theorem \ref{QpairforSegalspaces}.
In \ref{theoremA'}, we shall prove that

\medskip

\noindent{\bf Theorem} \emph{ The adjoint pair of functors 
$$t_!:{\bf S}^{(2)} \leftrightarrow {\bf S}:t^!$$
is a Quillen equivalence
between the model category for complete Segal spaces 
and the model category for quasi-categories.}

\medskip

We have stated the main results of the section.
We now proceed to the proofs.

\bigskip

Let $u_0$ be the inclusion $\{0\}\subset J$.

\begin{lemma}  \label{compSegspace} \cite{Rez} A Segal space
$X$ is complete iff
the map $$u_0\backslash X: J\backslash X  \longrightarrow 1 \backslash X$$
is a trivial fibration. 
\end{lemma} 

\noindent{\bf Proof}: 
By definition, a Segal space $X$
is complete  if 
the map
$$t\backslash X: 1\backslash X  \longrightarrow J \backslash X$$
obtained from the map $t:J\to 1$
is a weak homotopy equivalence.
But we have $(u_0\backslash X)(t\backslash X)=id$
since we have $tu_0=id$. Hence the map $t\backslash X$
is a weak homotopy equivalence iff the map $u_0\backslash X$
is a weak homotopy equivalence by three-for-two.
But the map $u_0\backslash X$ is a Kan fibration by 
\ref{Reedyfibrationvertically}. Hence the map $u_0\backslash X$
is a weak homotopy equivalence iff it is a trivial fibration.
\QED

\begin{lemma}\label{sliceofRezkfibrations} Let $f: X \ra Y$ be a Rezk fibration 
between complete Segal spaces.
then the map
$$\langle f/v \rangle :X/T\longrightarrow Y/T\times_{Y/S}X/S$$
is a quasi-fibration for any monomorphism
of simplicial sets $v:S\to T$. 
\end{lemma}

\noindent
{\bf Proof}: The map $\langle f/v \rangle$ is a mid fibration
between quasi-categories by \ref{v-fibSegalspace}.
Hence it suffices to show
that it has the right
lifting property with respect
to the inclusion $u_0:\{0\}\subset J$ by \ref{qfbetqc}.
But the condition
$u_0\pitchfork \langle f/v\rangle$
is equivalent to the condition
$v\pitchfork \langle u_0\backslash f\rangle$
by \ref{lift2}.
Hence it suffices to show
that the map
$$\langle u_0\backslash f\rangle
:J\backslash X \longrightarrow J\backslash Y \times_{1 \backslash Y  }1\backslash X $$
is a trivial fibration.
But it
 is a Kan fibration by \ref{Reedyfibrationvertically}.
 Hence it suffices to show
 that it is a weak homotopy equivalence.
But the horizontal maps in the commutative square 
 $$
\xymatrix{
J \backslash X  \ar[d] \ar[r]& 1\backslash X  \ar[d]\\
 J \backslash Y \ar[r]  & 1 \backslash Y     }
$$
are
trivial fibrations by \ref{compSegspace}.
It follows that 
$\langle u_0\backslash f\rangle$
is a weak homotopy equivalence.
It is thus a trivial fibration.
We have proved that $\langle f/v \rangle$ 
is a quasi-fibration.
\QED

\begin{proposition} \label{CompleteSegalspace2} A bisimplicial set $X$
is a complete Segal space
iff the following two conditions are satisfied:
\begin{itemize}
\item {\rm (i)} the map
$X\slash \delta_n $
is a quasi-fibration for every $n\geq 0$;
\item {\rm (ii)} $X$ is categorically constant.
\end{itemize}
\end{proposition}

\noindent{\bf Proof}: ($\Rightarrow$) A complete Segal space $X$  
is categorically constant by \ref{Segalspace3}.
Moreover, the map $X\slash \delta_n $
is a quasi-fibration for every $n\geq 0$ by \ref{sliceofRezkfibrations},
since the map $X\to 1$ is a Rezk fibration.
 ($\Leftarrow$) The bisimplicial set $X$ is h-fibrant by condition (i). 
 Let us show
that it is v-fibrant. 
By \ref{Reedyfibrationvertically} it suffices to show 
that the map $X/v$ is a trivial fibration for every
anodyne map $v\in {\bf S}$.
Let ${\cal A}$ be the class of monomorphisms
$v\in {\bf S}$ such that the map $X/v$ is a trivial fibration.
Let us show that ${\cal A}$ is saturated.
The condition $\delta_n\pitchfork (X\slash v)$
is equivalent to the condition $v\pitchfork (\delta_n\backslash X)$
by \ref{lift2}. Thus, a monomorphism $v$ belongs to ${\cal A}$ iff we have $v\pitchfork (\delta_n\backslash X)$
for every $n\geq 0$. It follows that the class ${\cal A}$ is saturated.
Let us show that ${\cal A}$ has the right cancellation property.
The map $X\slash v$
is a quasi-fibration for any monomorphism $v$ by \ref{Reedyfib1}
since $X$ is h-fibrant.
Thus, $X\slash v$ is a trivial fibration iff it is a weak categorical equivalence.
It follows by three-for-two that ${\cal A}$ has the right cancellation property.  
Let us show that every face map $d_i:\Delta[n-1]\to \Delta[n]$ belongs to  ${\cal A}$.
The canonical maps $X_{\star 0}\to X_{\star n}$ and $X_{\star 0}\to X_{\star, n-1}$ are weak categorical equivalences, since $X$
is categorically constant. It follows by three-for-two that the map $X/d_i:X_{\star n}\to X_{\star,n-1}$
is a weak categorical equivalence. This proves that $d_i\in {\cal A}$.
It then follows by Lemma \ref{Anodyne} 
that every anodyne map belongs to $\cal A$.
This shows that $X$ is v-fibrant.
Thus, $X_0$ is a Kan complex by \ref{Reedyfact2}.
It follows that $X$ is a Segal space by \ref{Segalspace2}.
Let us show that the Segal space $X$ is complete. 
The inclusion $u_0:\{0\}\subset J$ is an equivalence of categories.
It is thus a weak categorical equivalence.
Hence the map $u_0\backslash X: J\backslash X  \longrightarrow 1 \backslash X$
is a trivial fibration by \ref{Reedyfib1} since $X$ is h-fibrant.
This shows that $X$
is a complete Segal space by \ref{compSegspace}.
\QED

\begin{theorem} \label{Rowise2}The Rezk model structure $({\cal C}_R, {\cal W}_R,{\cal F}_R)$ is a Bousfield localisation of
the horizontal model structure $({\cal C}^h_1, {\cal W}^h_1,{\cal F}^h_1)$.
An h-fibrant bisimplicial set  is a complete Segal space
iff it is categorically constant. 
A row-wise weak categorical equivalence is a 
Rezk weak equivalence.
\end{theorem}

For the notion of Bousfield localisation, see Definition \ref{DefBousfieldLocalisationA}. 

\medskip

\noindent{\bf Proof}: Let us prove the first statement.
We have ${\cal C}^h_1={\cal C}_R$,
since  ${\cal C}^h_1$ is the class of monomorphisms by \ref{ReedyQuasi-categoriesmodel}.
 If $f:X\to Y$
is a Rezk fibration between complete Segal space,
then the map $\langle f/\delta_n\rangle\in  {\bf S}$ is a quasi-fibration
for every $n\geq 0$ by \ref{sliceofRezkfibrations}.
This means 
that a Rezk fibration between complete Segal space 
is an h-fibration.
This proves the first statement by \ref{Quillenpair2}. 
The second statement follows from \ref{CompleteSegalspace2}.
The third statement is a consequence
of the first and of Proposition \ref{KBlemma}.
\QED

In particular, a map between complete Segal spaces
is a Rezk weak equivalence (resp. a Rezk fibration)  iff it is a row-wise weak categorical equivalence
(resp. an h-fibration).

\begin{proposition}\label{BoxleftQuillen} The box product functor
$\Box:{\bf S}\times {\bf S}\to {\bf S}^{(2)}$
is a left Quillen functor of two variables
$$({\cal C}_1, {\cal W}_1,{\cal F}_1)\times ({\cal C}_0, {\cal W}_0,{\cal F}_0)\to
({\cal C}_R, {\cal W}_R,{\cal F}_R),$$
where $({\cal C}_1, {\cal W}_1,{\cal F}_1)$ is the model structure for quasi-categories  on $ {\bf S}$, $({\cal C}_0, {\cal W}_0,{\cal F}_0)$
is the classical model structure  on $ {\bf S}$
and $({\cal C}_R, {\cal W}_R,{\cal F}_R)$ is the Rezk model structure on $ {\bf S}^{(2)}$.
\end{proposition}

\noindent
{\bf Proof}: 
 Observe that the cofibrations are the monomorphisms
 in the three model structures,
Let $u:A\to B$ and $v:S\to T$ be a pair of monomorphisms in $\bf S$.
The map $u\Box' v$ is a cofibration since it is monic.
If $v\in {\cal W}_0$, let us show that $u\Box' v\in  {\cal W}_R$.
But we have $u\Box' v\in {\cal W}'_0$
by \ref{Reedycofibration2}. 
The result follows 
since a column-wise weak homotopy equivalence is a Rezk weak equivalence.
If $u\in {\cal W}_1$, let us show that $u\Box' v\in  {\cal W}_R$.
For this it suffices to show that we have $(u\Box' v)\pitchfork f$
for every Rezk fibration between complete Segal spaces $f:X\to Y$ by
\ref{lift3}. 
But the condition
$(u\Box' v)\pitchfork f$
is equivalent to the condition
$u\pitchfork \langle f/v\rangle$
by \ref{lift2}.
The map $\langle f/v\rangle$ is a quasi-fibration by \ref{Reedyfib1}.
Hence we have $u\pitchfork \langle f/v\rangle$,
since $u\in  {\cal C }_1\cap {\cal W}_1$.
\QED

We recall that the first projection $p_1: \Delta\times \Delta \to \Delta$
is left adjoint to the functor $i_1:\Delta \to  \Delta\times \Delta$ 
defined by putting 
$i_1([n])=([n],[0])$ for every $n\geq 0$.

\begin{proposition} \label{leftQ1}
The adjoint pair of functors $$p^*_1: {\bf S}\rightarrow  {\bf S}^{(2)}:i_1^*$$
is a homotopy localisation
between the model category for quasi-categories
and the model category for complete Segal spaces.
\end{proposition}

See Definition \ref{Localisation1} for the notion of homotopy localisation.
\medskip

\noindent {\bf Proof}: 
We have $p_1^*(A)=A\Box 1$ for every $A\in {\bf S}$.
The functor $A\mapsto A\Box 1$ 
is a left Quillen functor since the box product functor $(A,B)\mapsto A\Box B$ is a left Quillen functor 
by \ref{BoxleftQuillen} and since $1$ is cofibrant.
This proves that the pair $(p_1^*,i_1^*)$ is a Quillen pair.
Let us show that it is a homotopy localisation. For this, we shall use
Proposition \ref{Localisation2}.
It suffices to show that 
the adjunction counit $\epsilon:p^*_1i^*_1X\to X$
is a Rezk weak equivalence for every complete Segal space $X$.
For this, it suffices to show by Theorem \ref{Rowise2} that 
the map $\epsilon_{\star n}:X_{\star 0}\to X_{\star n}$ 
 is a weak categorical equivalence for every $n\geq 0$.
But $\epsilon_n$ is equal to the map $X_{\star 0}\to X_{\star n}$
 obtained from the map $[n]\to [0]$.
This proves the result since $X$ is categorically
 constant by \ref{CompleteSegalspace2}.
 \QED

\bigskip

If $X$ is a bisimplicial set, then for any pair of simplicial sets $A$ and $B$
 there is a natural bijection 
between the maps $A\to X/B$
and the maps $B\to A\backslash X$.
This means that 
 the contravariant functors 
$$A\mapsto A\backslash X\quad {\rm and}\quad B\mapsto X/B$$ 
(from $\bf S$ to itself) are mutually right adjoint.

\begin{lemma}  \label{continuousfunct}
Every continuous functor $G:{\bf S}^o\to {\bf S}$
is of the form $G(A)=A \backslash X$
for a simplicial space $X\in {\bf S}^{(2)}$.
We have
$X_{m}=G(\Delta[m])$
for every $m\geq 0$.
\end{lemma}

\noindent{\bf Proof}:  Let $X$
be the simplicial space defined by putting 
$X_{m}=G(\Delta[m])$
for every $m\geq 0$.
Then we have 
$\Delta[m]\backslash X=X_m=G(\Delta[m])$
for every $m\geq 0$. 
It follows that we have $A\backslash X=G(A)$
for every simplicial set $A$,
since the functors $A\mapsto A\backslash X$
and $A\mapsto G(A)$ are both continuous
and every simplicial set is a colimit 
of simplices.
\QED

If $X$ is a quasi-category, we shall denote by 
$\Gamma(X)$
the simplicial space obtained by putting
$$\Gamma(X)_{m}=J(X^{\Delta[m]})$$
for every $m\geq 0$.

\begin{proposition} \label{Double1} 
There are canonical isomorphisms
$$A\backslash \Gamma(X)  =J(X^A)\quad {\rm and } \quad \Gamma(X)/A=X^{(A)},$$
natural in $A\in {\bf S}$.
\end{proposition}

\noindent {\bf Proof}: 
The contravariant functors $A\mapsto J(X^A)$
and $A\mapsto X^{(A)}$ are mutually right adjoint by \ref{Jfunctoradjoint}.
Hence the contravariant functor $A\mapsto J(X^A)$
is continuous. It follows by \ref{continuousfunct}
that 
we have $A\backslash \Gamma(X)  =J(X^A)$
for every simplicial set $A$.
But the contravariant functors $A\mapsto A\backslash \Gamma(X)$
and $A\mapsto \Gamma(X)/A$ are mutually right adjoint.
Hence we have  $\Gamma(X)/A=X^{(A)}$ for every simplicial set $A$
by the uniqueness of a right adjoint.
\QED

\medskip

If $X$ is a simplicial set, then we have
$$i_1^*\Gamma(X)=\Gamma(X)_{\star 0}=\Gamma(X)/1=X^{(1)}=X$$
by \ref{Double1}.
By the adjointness $p_1^*\dashv i_1^*$
we obtain a natural map
$p_1^*(X)\to  \Gamma(X)$.

\begin{proposition}\label{fibrantreplacementResz}  If $X$ is a quasi-category,
then $\Gamma(X)$ is a complete Segal space and 
the natural map $p_1^*(X)\to  \Gamma(X)$
is a Rezk weak equivalence.
Hence $\Gamma(X)$
is a fibrant replacement of the
simplicial space $p_1^*(X)=X\Box 1$.
\end{proposition}

\noindent {\bf Proof}: Let $X$
be a quasi-category. Let us first show that $\Gamma(X)$ is vertically Reedy fibrant.
Let $\delta_n$ be the inclusion $\partial \Delta[n] \subset\Delta[n ]$.
The map $X^{\delta_n}:X^{\Delta[n]}\to X^{\partial \Delta[n]}$ is a quasi-fibration 
since the model structure for quasi-categories is 
cartesian closed by \ref{QCatmodel}.
Hence the map $J(X^{\delta_n}):J(X^{\Delta[n]})\to J(X^{\partial \Delta[n]})$
is a Kan fibration by \ref{Jfunctoradjoint}.
But $J(X^{\delta_n})$
is isomorphic to the map
$\delta_n\backslash \Gamma(X):\Delta[n]\backslash \Gamma(X)\to \partial \Delta[n]\backslash \Gamma(X)$
by \ref{Double1}. This shows that  the map
$\delta_n\backslash \Gamma(X)$
is a Kan fibration.
We have proved that  $\Gamma(X)$ is  vertically Reedy fibrant.
Let us now show that $\Gamma(X)$ is a Segal space.
The inclusion $i_n:I_n\subseteq \Delta[n ]$ is mid anodyne by \ref{midado2}.
It is thus a weak categorical equivalence by \ref{midcofibration}.
Hence the map
$X^{i_n}$ is a trivial fibration by \ref{QCatmodel}.
It follows that the map $J(X^{i_n})$
is a trivial fibration by \ref{Jfunctoradjoint}.
But $J(X^{i_n})$
is isomorphic to the map
$i_n\backslash \Gamma(X)$
by \ref{Double1}. 
This shows that  the map
$i_n\backslash \Gamma(X)$
is a trivial fibration.
We have proved that  $\Gamma(X)$ is a Segal space.
It remains to show that $\Gamma(X)$ is a
complete Segal space. 
The map $p:J\to 1$
is an equivalence of categories.
Hence the map $X^p$
is an equivalence of quasi-categories by \ref{QCatmodel}.
It follows that the map $J(X^p)$
is a homotopy equivalence by \ref{Jfunctoradjoint}.
But $J(X^p)$
is isomorphic to the map
$p\backslash \Gamma(X)$
by \ref{Double1}. 
This shows that the map 
$p\backslash \Gamma(X)$
is a homotopy equivalence.
We have proved that $\Gamma(X)$ is a
complete Segal space.
Let us now show that
the natural map $p_1^*(X)\to  \Gamma(X)$
is a Rezk weak equivalence.
For this, it suffices to show that the natural map $p_1^*(X)\to  \Gamma(X)$ 
is a row-wise weak categorical equivalence by Lemma \ref{Rowise2}. 
But the map $p_1^*(X)_{\star n}\to \Gamma(X)_{\star n}$
is equal to the map $X^{(t_n)}:X^{(1)} \to  X^{(\Delta[n])}$,
where $t_n:\Delta[n]\to 1$.
But $X^{(t_n)}$  is an equivalence of quasi-categories
by \ref{Jdual2}, since $t_n$
is a weak homotopy equivalence.
\QED
 
\bigskip

\begin{theorem}  \label{Theorem A}
The adjoint pair of functors $$p^*_1: {\bf S}\rightarrow  {\bf S}^{(2)}:i_1^*$$
is a Quillen equivalence
between the model category for quasi-categories
and the model category for complete Segal spaces.
\end{theorem}

\noindent{\bf Proof}: We shall use propostion \ref{Quillenequiv}.
We saw in \ref{leftQ1} that the pair $(p^*_1, i_1^*)$
is a homotopy  localisation.
Hence it suffices to show that the pair $(p^*_1, i_1^*)$
is a homotopy colocalisation by \ref{Quillenequiv}. 
By \ref{Localisation2}, for this it suffices to show that 
the map $A \ra i_1^*Rp^*_1A$ is a weak
equivalence for every fibrant-cofibrant object $A \in {\bf S}$, where $p^*_1A\to Rp^*_1A$ denotes a fibrant replacement of $p^*_1A$. But
we can take $Rp^*_1A =\Gamma(A)$
by \ref{fibrantreplacementResz}.  In this case we have
$$i^*_1Rp^*_1A=i_1^*\Gamma(A)=A$$
and the canonical map $A\to i^*_1Rp^*_1A $ is
the identity. The result is proved.
\QED

\bigskip

Recall the ``total space" functor 
$t_!:{\bf S}^{(2)} \to {\bf S}$
of Theorem \ref{QpairforSegalspaces}.

\begin{theorem}  \label{theoremA'} 
The adjoint pair of functors 
$$t_!:{\bf S}^{(2)} \leftrightarrow {\bf S}:t^!$$
is a Quillen equivalence
between the model category for complete Segal spaces 
and the model category for quasi-categories.
\end{theorem}

\noindent{\bf Proof}: Let us first show that $(t_!,t^!)$
is a Quillen pair. We saw in \ref{QpairforSegalspaces} 
that it is a Quillen pair
between the model category for Segal spaces 
and the model category for quasi-categories.
Hence it suffices to show by \ref{Quillenpair2}.
that
the functor $t^!$ takes a quasi-category to a complete Segal space.
If $X$ is a quasi-category then $t^!(X)$ 
is a Segal space by \ref{Qpairforhorizontalandvertical}.
Let us show that the Segal space $t^!(X)$ 
is complete.  For this, it suffices to show that the map 
$u_0\backslash t^!(X)$
is a trivial fibration by \ref{compSegspace},
where $u_0$ denotes the inclusion $\{0\}\subset J$.
But the map $u_0\backslash t^!(X)$ 
is isomorphic to the map
$k^!(X^{u_0})$
by Lemma \ref{Boxtoproductandt}. 
Hence it suffices to show that the map $k^!(X^{u_0})$
is a trivial fibration.
But $u_0$ is a weak categorical equivalence
since it is an equivalence of categories.
It follows that the map $X^{u_0}$
is a trivial fibration by \ref{QCatmodel}
since $u_0$ is monic.
This shows that the map 
$u_0\backslash t^!(X)$
is a trivial fibration.
We have proved that $t^!(X)$ is 
a complete Segal space. 
We have proved that $(t_!,t^!)$
is a Quillen pair
between the model category for complete Segal spaces 
and the model category for quasi-categories.
It remains to show that
it is a Quillen equivalence.
The composite $t_!p_1^*: {\bf S}\to  {\bf S}$
is isomorphic to the identity functor
since we have
$$t_!p_1^*(A)= t_!(A\Box 1)=A\times k_!(1)= A$$
for every simplicial set $A$ by \ref{Boxtoproductandt}
and since $k_!(1)=1$.
Hence the composite 
$i_1^*t^!: {\bf S}\to  {\bf S}$
is also isomorphic to the identity functor
by adjointeness. 
We saw in \ref{Theorem A} that the pair $(p_1^*,i_1^*)$
is a Quillen equivalence.
It follows by three-for-two in \ref{Quillenequiv3for2}
that the pair $(t_!,t^!)$
 is a Quillen equivalence.\QED

\section {Two equivalences with Segal categories }

Segal categories were first introduced by Dwyer, Kan
and Smith \cite{DKS}, where they are called special $\Delta^o$-diagrams
of simplicial sets .
The theory of Segal categories was extensively developed by 
Hirschowitz and Simpson for application to algebraic geometry.
A simplicial space $X$ is called a precategory
if its  first column $X_0$ is discrete. There is a model structure
on the category of precategories
in which the fibrant objects are Segal categories.
The goal of the section is to
establish two Quillen equivalences 
$$q^*:{\bf S}\leftrightarrow {\bf PCat}:j^* \quad {\rm and} \quad d^*:{\bf PCat}\leftrightarrow  {\bf S}:d_*$$ 
between the model category for quasi-categories
and the model category for Segal categories.
The functor $j^*$
associates to a precategory $X$ its first row $X_{\star 0}$.
The functor $d^*$ associates to a precategory its diagonal.

\bigskip

We recall that a simplicial space $X:\Delta^o\to {\bf S}$
is called a {\it precategory} if 
$X_0$ is discrete. 
We shall denote by ${\bf PCat}$
the full subcategory
of ${\bf S}^{(2)}$
spanned by the precategories.

\medskip

Consider the functor $i_2:\Delta \to \Delta\times \Delta$ defined by putting $i_2([n])=([0],[n])$.
The functor $i_2^*: {\bf S}^{(2)}\to {\bf S}$ associates to a simplicial
space $X$ its first column $X_0$. 
The functor $i_2^*$ has a right adjoint $Cosk={(i_2)}_\star:  {\bf S}\to {\bf S}^{(2)}$.
If $A$ is a simplicial set, then we have
$${ Cosk}(A)_n=A^{[n]_0}$$
for every $n\geq 0$,
where $[n]_0$ denotes the set of vertices on $\Delta[n]$.
If $X$ is a simplicial space, the unit of the adjunction $i_2^*\dashv {(i_2)}_\star$
is a canonical map 
$$v_X:X\to{Cosk}(X_0) $$ called the {\it vertex map}.
Let us suppose that $X$ is a precategory.
Then the map $(v_X)_n:X_n\to{Cosk}(X_0)_n$
takes its values in a discrete simplicial set $ X_0^{[n]_0}$
for each $n\geq 0$. We thus have a decomposition 
$$ X_n=\bigsqcup_{a\in X_0^{[n]_0}} X(a),$$
where $X(a)=X(a_0,a_1,\ldots,a_n)$ denotes the fiber of the vertex map
$(v_X)_n: X_n \to X_0^{[n]_0}$
at $a=(a_0,a_1,\cdots,a_n)$.
If $u:[m]\to [n]$ is a map in $\Delta$, then the map $X(u):X_n\to X_m$
induces a map
$$X(a_0,a_1,\ldots,a_n)\to X(a_{u(0)},a_{u(1)},\ldots,a_{u(m)})$$
for every $a\in  X_0^{[n]_0}$. 

\medskip

A precategory $X$ is called a {\it Segal category} in \cite{HS} if it satisfies the Segal condition
\ref{defSegalcondition}. It is easy to verify that we have a decomposition
$$I_n\backslash X =\bigsqcup_{a\in X_0^{[n]_0}} X(a_0,a_1)\times \cdots \times X(a_{n-1},a_n).$$
It follows that a precategory $X$ is a Segal category iff
the canonical map
$$X(a_0,a_1,\ldots,a_n)\to X(a_0,a_1)\times \cdots \times X(a_{n-1},a_n)$$
is a weak homotopy equivalence for every $n\geq 2$
and every $a\in  X_0^{[n]_0}$ (the condition is trivially satisfied
if $n<2$).

\medskip

\begin{definition} \label{deffullyfaithful}  \cite{HS} A map of Segal categories
$f:X\to Y$ is said to be
{\it fully faithful} if 
the map 
$$X(a,b)\to Y(fa,fb)$$
is a weak homotopy equivalence for every pair $a,b\in X_0.$
\end{definition}

If $C$ is a small category, then the bisimplicial set $N(C)=C\Box 1$
is a Segal category.
The functor $N:{\bf Cat}\to {\bf PCat}$
has a left adjoint
$$\tau_1:{\bf PCat} \to  {\bf Cat}.$$
We say that
$\tau_1 X$ is the {\it fundamental category} of a precategory $X$.
We shall say that a map of precategories
$f:X\to Y$
is {\it essentially surjective} if the functor $\tau_1(f):\tau_1X\to \tau_1 Y$
is essentially surjective.

\begin{definition} \label{defessensurj}  \cite{HS} A map between Segal categories
$f:X\to Y$ is called an {\it equivalence} if it is  fully faithful
and essentially surjective.
\end{definition}

\medskip

Hirschowitz and Simpson construct a completion functor $Seg: {PCat}\to {PCat}$
which associates to a precategory $X$ a Segal category $Seg(X)$ ``generated" by $X$.
A map of precategories $f:X\to Y$ is called a {\it weak categorical equivalence}
if the map $Seg(f):Seg(X)\to Seg(Y)$ is an equivalence of Segal categories.

\bigskip

\begin{theorem} [Hirschowitz-Simpson, see \cite{HS} ]  \label{HSmodel} The category of precategories 
${\bf PCat}$ admits a model structure in which a cofibration
is a monomorphism and a weak equivalence is a weak categorical equivalence. 
The model structure is left proper and cartesian closed.
\end{theorem}

We call  the model structure, the  {\it Hirschowitz-Simpson model structure}
or the {\it model structure for Segal categories}.
The model structure is cartesian closed by a result of
Pellisier in \cite{P}. 
A precategory is fibrant iff it is a Segal space by \cite{B3} and \cite{J4}.

\bigskip

A simplicial set $X: \Delta^o\to {\bf Set} $ is discrete
iff it takes every map in $ \Delta$ to a bijection.
It follows that a bisimplicial set $X: (\Delta^2)^o \to {\bf Set} $
is a precategory iff it takes every 
map in $[0]\times \Delta$ to a bijection.
Let us put 
$$\Delta^{\mid 2} =([0]\times \Delta)^{-1}(\Delta\times \Delta)$$ 
and let $\pi$ be the canonical functor
$\Delta^2\to \Delta^{\mid 2}$.

\begin{proposition} \label{pilowerstarupperatar}The functor $\pi^*$
induces an isomorphism between the presheaf category $[\Delta^{\mid 2},{\bf Set}] $ 
and the subcategory ${\bf PCat}\subset {\bf S}^{(2)}$.
\end{proposition}

We shall regard the functor $\pi^*$ as an inclusion
by adopting the same notation
for a contravariant functor $X:\Delta^{\mid 2}\to {\bf Set}$ and the precategory $\pi^*(X)$.
The functor $\pi^*: {\bf PCat}\subset {\bf S}^{(2)}$
has a left adjoint $\pi_!$
and a right adjoint $\pi_*$.

\medskip

\begin{theorem} [Bergner, see\cite{B2}] \label{Qequiv3} The pair of adjoint
functors $$\pi^*:{\bf PCat}\leftrightarrow  {\bf S}^{(2)}:\pi_*$$ 
is a Quillen equivalence between the model structure for Segal categories
and  the Resz model structure on $ {\bf S}^{(2)}$. 
A map of precategories $X\to Y$
is a weak categorical equivalence
iff the map $\pi^*(u)$ is a Rezk weak equivalence.
\end{theorem}

\bigskip

The functor $i_i:\Delta\to \Delta \times \Delta$ defined by putting $i_1([n])=([n],0)$
is right adjoint to the projection $p_1:\Delta \times \Delta \to \Delta$.
The projection $p_1$
inverts every arrow in $[0]\times \Delta$.
Hence there is a unique functor $q:\Delta^{\mid 2}\to \Delta$
such that $q\pi=p_1$.
The functor $j=\pi i_1:\Delta \to \Delta^{\mid 2}$
is then right adjoint to the functor $q$.
Hence the functor
$j^*:{\bf PCat}\to{\bf S}$ 
is right adjoint to the functor
$q^*$.  If $X$ is a precategory, then $j^*(X)$
 is the first row of $X$. 
 If $A\in {\bf S}$, then $q^*(A)=A \Box 1$.
 The following result was conjectured by Bertrand T\"oen in \cite{T1}:

\bigskip

\begin{theorem} \label{ConjectureToen}
The adjoint pair of functors 
$$q^*:{\bf S}\leftrightarrow {\bf PCat}:j^*$$ 
is a Quillen equivalence
between the model category for quasi-categories
and the model category for Segal categories.
\end{theorem} 

\bigskip

\noindent{\bf Proof}: Let us show that $q^*$
is a left Quillen functor. Obviously, $q^*$
preserves monomorphisms .
Let us show that
it takes a weak categorical equivalence
to a weak categorical equivalence.
The functor $p_1^*$ 
takes a weak categorical equivalence
to a Rezk weak equivalence by \ref{Theorem A}
(and by Proposition \ref{KBlemma}).
Thus, if $u\in {\bf S}$ is a weak categorical equivalence,
then $p_1^*(u)$ is a Rezk weak equivalence.
But we have $p_1^*= \pi^*q^*$,
since we have $q\pi=p_1$.
Thus, $\pi^*q^*(u)$ 
is a Rezk weak equivalence.
It follows by \ref{Qequiv3} that $q^*(u)$ 
is a weak categorical equivalence. 
We have proved that $q^*$ is a left Quillen functor.
It remains to show that the pair $(q^*,j^*)$
it is a Quillen equivalence.
But the pair $(p_1^*,i_1^*)$
is a Quillen equivalence by \ref{Theorem A}. 
Hence also the pair $(q^*,j^*)$ by three-for-two in \ref{Quillenequiv3for2},
since the pair $(\pi^*,\pi_*)$ is a Quillen equivalence by \ref{Qequiv3}
and since we have $\pi^*q^*=p_1^*$. \QED

\medskip

Let us put $d=\pi\delta:\Delta \to \Delta^{\mid 2}$,
where $\delta$ is the diagonal functor $\Delta \to \Delta \times \Delta$.
The simplicial set $d^*(X)$
 is the diagonal of a precategory $X$. 
The functor $$d^*: {\bf PCat}\to  {\bf S}$$
admits a left adjoint $d_!$ and a right adjoint $d_*$.

\bigskip

\begin{theorem} \label{diagonalequivforprecategories}
The adjoint pair of functors 
$$d^*:{\bf PCat}\leftrightarrow {\bf S} :d_*$$ 
is a Quillen equivalence
between the model category for Segal categories
and the model category for quasi-categories.
\end{theorem}

The proof is given after Lemma \ref{gammafunctor and qf}.
Let  $\delta_\star $ be the right adjoint of the 
functor $ \delta^*: {\bf S}\to  {\bf S}^{(2)}$.

\begin{lemma}  \label{diagonalbox}  For every  $A,B,X\in {\bf S}$ we have
$$\delta^*(A\Box B)=A\times B, \quad \quad  A\backslash \delta_*(X) =X^A
 \quad{\rm and}\quad\delta_*(X)/B =X^B.$$
\end{lemma}  

\noindent{\bf Proof}: We have 
$$\delta^*(A\Box 1)=\delta^*p_1^*(A)=(p_1\delta)^*( A)=A$$
since $p_1\delta =id$. Similarly, we have $\delta^*(1\Box B)=B$
since $p_2\delta =id$.
But $$A\Box B=(A\Box 1)\times (1\Box B).$$
Thus,
$\delta^*(A\Box B)=\delta^*(A\Box1)\times \delta^*(1\Box B)=A\times B,$
since the functor $\delta^*$ preserves products.
Let us show that $A\backslash \delta_*(X) =X^A$.
If $B$ is a simplicial set, there is a natural bijection
between the maps $B\to A\backslash \delta_*(X) $,
the maps $A\Box B\to \delta^*(X) $, the maps $\delta^*(A\Box B)\to X$, 
 the maps $A\times B\to X$ and the maps $B\to X^A$.
 This shows by Yoneda lemma that $A\backslash \delta_*(X) =X^A$.
The formula $\delta_*(X)/B =X^B$
is proved similarly.
\QED

\begin{proposition} \label{pilowerstar} If $X$ is a simplicial space, then we have 
a pullback square of bisimplicial sets,
$$
\xymatrix{
\pi_*X  \ar[d]\ar[r] & X \ar[d]^v\\
{ Cosk}(X_{00}) \ar[r] & { Cosk}(X_0),}
$$
where $v=v_X:X\to Cosk(X_0)$ is the vertex map.
\end{proposition}

\noindent{\bf Proof}: 
Let $P$ be the bisimplicial set
defined by the pullback square
$$
\xymatrix{
P  \ar[d]\ar[r]^i & X \ar[d]^v\\
{ Cosk}(X_{00}) \ar[r] & { Cosk}(X_0).}
$$
The map $v$ induces an isomorphism 
on the first columns, hence also the map
$P\to { Cosk}(X_{00})$.  Thus, $P_0=X_{00}$.
This shows that $P$ is a precategory.
Let us show that the map $i:P\to X$ coreflects $X$
in the subcategory  ${\bf PCat}$.
For this, we have to show that if 
$Z$ is a precategory, then
every map $f:Z\to X$ factors uniquely through $i$.
But the map
$f_0:Z_0\to X_0$ factors through the inclusion  $X_{00}\subseteq X_0$
since $Z_0$ is discrete.
The result then follows by using the adjunction $(i_2)^*\dashv (i_2)_\star$.
\QED

Recall from \ref{fibrantreplacementResz} the functor $\Gamma:{\bf QCat}\to {\bf S}^{(2)}$
which associates to a quasi-category $X$ a complete Segal space $\Gamma(X)$.
We have
$$\Gamma(X)_{\star n}=\Gamma(X)/\Delta[n]=X^{(\Delta[n])}$$
for every $n\geq 0$ by \ref{Double1}.

\begin{lemma}  \label{rightbergnerfunctor2} If $X$ is a quasi-category, then 
$d_*X=\pi_*\Gamma (X).$
\end{lemma}  

\noindent{\bf Proof}: 
It suffices to show
that we have $(\pi_*\Gamma X)_n= (\pi_*\delta_*X)_n$
for every $n\geq 0$.
If $Y$ is a simplicial space,
then we have a pullback square of simplicial sets
$$
\xymatrix{
(\pi_*Y)_n  \ar[d]\ar[r] & Y_n  \ar[d]\\
Y_{00}^{[n]_0} \ar[r] & Y_0^{[n]_0}}
$$
by Lemma \ref{pilowerstar}.
In particular, if $Y=\delta_*X$, we have a 
pullback square of simplicial sets
$$
\xymatrix{
(\pi_*\delta_*X)_n \ar[d]\ar[r] & X^{\Delta[n]} \ar[d]\\
X_{0}^{[n]_0} \ar[r] & X^{[n]_0}. }
$$
since $Y_n=\Delta[n]\backslash \delta_*X=X^{\Delta[n]}$
by \ref{diagonalbox}.
The projection $ X^{\Delta[n]} \to X^{[n]_0}$
is a conservative quasi-fibration by \ref{conservativeprojection}.
Hence also the projection $(\pi_*\delta_*X)_n \to X_{0}^{[n]_0}$
by \ref{basechangeconservativeJ}. It follows that  the simplicial set $(\pi_*\delta_*X)_n$
is a Kan complex by \ref{conservativefibKanfibration},
since $X_{0}^{[n]_0} $ is a Kan complex.
If we apply the functor $J$ to the pullback square above
we obtain a pullback square
$$
\xymatrix{
(\pi_*\delta_*X)_n \ar[d]\ar[r] & J(X^{\Delta[n]})  \ar[d] \\
X_{0}^{[n]_0} \ar[r] & J(X)^{[n]_0}, }
$$
since the functor $J$ preserves pullbacks by \ref{Jfunctoradjoint}   
and since the vertical map on the left hand side
is a map between Kan complexes. 
But we have 
$J(X^{\Delta[n]})= (\Gamma X)_n$, $J(X)=(\Gamma X)_0$
and $X_{0}=(\Gamma X)_{00}$ by \ref{Double1}. 
We thus obtain a pullback square
$$
\xymatrix{
(\pi_*\delta_*X)_n \ar[d]\ar[r] & (\Gamma X)_n  \ar[d]\\
(\Gamma X)_{00}^{[n]_0}\ar[r] & (\Gamma X)^{[n]_0} }.
$$
This shows that 
$(\pi_*\delta_*X)_n = (\pi_*\Gamma X)_n$.
\QED

\begin{lemma}  \label{gammafunctor and qf} If $f:X\to Y$ is a quasi-fibration between
quasi-categories, then the map $\Gamma(f): \Gamma X \to \Gamma Y$
is a Rezk fibration between complete Segal spaces.
\end{lemma}

\noindent{\bf Proof}: The bisimplicial sets
$\Gamma X$ and $\Gamma Y$
are complete Segal spaces by \ref{fibrantreplacementResz}. 
Hence it suffices to show that $\Gamma(f)$
is a v-fibration. For this it suffices to show that
the map $u\backslash \Gamma (f)$
is a Kan fibration for every monomorphism $u:A\to B$.
But the square
$$\xymatrix{
B\backslash \Gamma X \ar[r]\ar[d]& A\backslash \Gamma X \ar[d] \\
B\backslash\Gamma Y   \ar[r]& A\backslash \Gamma Y
 }
 $$
is isomorphic to the square 
$$\xymatrix{
 J(X^B) \ar[r]\ar[d]& J(X^A ) \ar[d] \\
 J(Y^B)   \ar[r]&  J(Y^A)
 }
 $$
by Lemma \ref{Double1} .
Hence the map $ u\backslash \Gamma (f)$
is isomorphic to the image by the functor $J$
of the map
 $$\langle u,f\rangle:X^B\to Y^B\times_{Y^A}X^A$$
But $\langle u,f\rangle$ is a quasi-fibration
by Theorem \ref{QCatmodel}. 
It follows that $J\langle u,f\rangle$
is a Kan fibration by \ref{Jfunctoradjoint}.
We have proved that $\Gamma(f)$
is a Rezk fibration.\QED

\noindent{\bf Proof of Theorem  \ref{diagonalequivforprecategories}}: 
Let us first show that
the pair $(d^*,d_*)$ is a Quillen pair.
For this, we shall use
the criteria of \ref{Quillenpair2}.
Obviously, the functor $d^*$ 
takes a monomorphism to a monomorphism.
Let us show that its right adjoint $d_*$
takes a fibration between  fibrant objects
to a fibration.
If $f:X\to Y$ is a quasi-fibration
between quasi-categories, let us show
that the map $d_*(f):d_*(X)\to d_*(Y)$
is a fibration.
But we have $d_*(f)=\pi_*\Gamma(f)$
by \ref{rightbergnerfunctor2}.
The map $\Gamma(f)$ is a Rezk fibration
by \ref{gammafunctor and qf}. Thus,
$\pi_*\Gamma(f)$ is a fibration since $\pi_*$
is a right Quillen functor by \ref{Qequiv3}.
We have proved that the pair $(d^*,d_*)$
is a Quillen pair. It remains to show that
it is a Quillen equivalence.
The composite $d^*q^*: {\bf S}\to  {\bf S}$
is isomorphic to the identity functor
since $qd=id$. We saw in \ref{ConjectureToen} that the pair $(q^*,j^*)$
is a Quillen equivalence.
It follows by three-for-two in \ref{Quillenequiv3for2}
that the pair $(d^*,d_*)$
 is a Quillen equivalence.\QED

\section {Addendum}

The model structure for quasi-categories in Theorem \ref{QCatmodel}
is not simplicial.
However, it is Quillen  equivalent to the model
category for complete Segal spaces, which is simplicial by Theorem \ref{Rezkmodel}.
In their paper {\it Simplicial structures on model categories and functors} \cite{RSS}
Rezk, Schwede and Shipley study the problem of associating 
to a model category $ {\cal E}$ a Quillen equivalent
simplicial model category.
A simplicial object $X:\Delta^o \to   {\cal E}$ is said to be
{\it homotopically constant} if it takes every map
in $\Delta$ to a weak equivalence in ${\cal E}$.
For any model category ${\cal E}$ 
there is at most one model structure $M^c=({\cal C}^c,  {\cal W}^c, {\cal F}^c)$
on the category $[\Delta^o, {\cal E}]$ such that
\begin{itemize}
\item{} The model structure $M^c$ is a Bousfield localisation of
the Reedy model structure on $[\Delta^o, {\cal E}]$;
\item{} The fibrant objects are the homotopically constant Reedy fibrant objects.
\end{itemize}
For the notion of Bousfield localisation, see Definition \ref{DefBousfieldLocalisationA}. 
When it exists, $M^c$ is called the {\it canonical model structure} on $[\Delta^o, {\cal E}]$.
Under certain conditions, the model category $M^c$ 
is shown to be simplicial and the
adjoint pair 
$$p^*: {\cal E}  \leftrightarrow [\Delta^{o}, {\cal E}] :i^*$$ 
to be a Quillen equivalence,
where $p^*(A)=1\Box A$ and $j^*(X)=X_0$.

\begin{theorem} \label{Addendumtheorem1} The Rezk model structure $({\cal C}_R, {\cal W}_R,{\cal F}_R)$ on ${\bf S}^{(2)}=[\Delta^o, {\bf S}]$ is
the canonical model structure 
associated to the model
structure for quasi-categories
$({\cal C}_1,  {\cal W}_1, {\cal F}_1)$ on ${\bf S}$.
\end{theorem}

\noindent{\bf Proof}: The model structure $({\cal C}_R, {\cal W}_R,{\cal F}_R)$ 
is a Bousfield localisation of the horizontal model structure 
$({\cal C}_1^h,  {\cal W}_1^h, {\cal F}_1^h)$
by \ref{Rowise2}. Moreover, an h-fibrant simplicial space $X$ is a complete Segal space 
iff it is categorically constant.
\QED

\section {Appendix}

The goal of this appendix is to review the basic homotopical algebra needed in the paper
and to introduce some notation.

\bigskip

We shall denote by $Ob{\cal C}$
the class of objects of a category ${\cal C}$ and  by ${\cal C}(A,B)$
the set of arrows between two objects of ${\cal C}$.
If $F:{\cal A}\to  {\cal B}$ and $G: {\cal B}\to {\cal A}$
are two functors, we shall write $F\dashv G$, or write
$$F: {\cal A}\leftrightarrow {\cal B}:G,$$
to indicate that the functor $F$
 is left adjoint to the functor $G.$

\medskip

We shall denote by ${\bf Set}$ the category of sets
and by ${\bf Cat}$ the category of small categories.
If $A$ is a small category and ${\cal E}$ is a category (possibly large) 
we shall denote the category of functors $A\to {\cal E}$
by $[A, {\cal E}]$ or by ${\cal E}^A$.
Recall that a {\it presheaf} on a small category $A$  is a 
contravariant functors $A\to {\bf Set}$.
We shall denote by $\hat A$  
the category $[A^o, {\bf Set}]$ of presheaves on $A$.
We shall regard the Yoneda functor
$y:A\to \hat A$ as an inclusion by adopting
the same notation for an object $a\in A$
and the representable functor $y(a)=A(-,a)$.
For every functor $u:A\to {\cal E}$, we shall denote
by $u^!:{\cal E}  \to \hat A$ the functor  obtained by putting $u^!(X)(a)={\cal E}(u(a),X)$
for every object $X\in {\cal E}$ and every object $a\in A$.
The functor $u^!$ has a left adjoint $u_!$
when the category ${\cal E}$ is cocomplete.
The functor $u_!:  \hat A  \to{\cal E}$
is the left Kan extension of the functor $u$ along
the Yoneda functor $A\to \hat A$. 
If $B $ is a small category and $u:A\to B$,
we shall denote by $u^*:\hat B\to \hat A$
the functor obtained by putting $u^*(X)=Xu$
for every presheaf $X\in \hat B$.
The functor $u^*$
has a left adjoint denoted $u_!$
and a right adjoint denoted $u_*$.

\medskip

We denote by ${\Delta }$ the category whose objects are the finite non-empty ordinals
and whose arrows are the order preserving maps.
The ordinal $n+1$ is represented by the ordered set $[n]=\{0,\ldots,n\}$,
so that $Ob\Delta=\{[n]:n\geq 0\}$.
A {\it simplicial set} is a presheaf on $\Delta$.  
If $X$ is a simplicial set, the set $X([n])$
is denoted by $X_n$ for every $n\geq 0$.
We denote the category of simplicial sets $\hat \Delta$ by
${\bf S}$. Recall that the simplex $\Delta[n]$
is defined to  be the representable functor $\Delta(-,[n] )$.
We shall denote its boundary by $\partial\Delta[n]$.
A category enriched over ${\bf S}$ is called
a {\it simplicial category}.

\medskip

Let $c: {\bf Set}\to {\bf S}$ 
be the functor which associate to a
set $S$ the constant simplicial set $cS$
obtained by putting $(cS)_n=S$
for every $n\geq 0$.
The functor
$c$ is full and faithful and we shall regard it as an 
inclusion ${\bf Set}\subset {\bf S}$
by adopting the same notation for $S$ and $cS$.
The functor $c$ has a left adjoint  
$$\pi_0:{\bf S}\to  {\bf Set},$$
where $\pi_0(X)$ is the set of connected components
of a simplicial set $X$.

\medskip

If $u:A\to B$
and $f:X\to Y$
are two maps in a category $\cal E$,
we write $u\pitchfork f$
to indicate that  $f$
has the right lifting property
with respect to $u$.
If $S$ is an object of $\cal E$, we write $u\pitchfork S$
to indicate that the map ${\cal E}(u,S):{\cal E}(B,S)\to {\cal E}(A,S)$
is surjective and we write  $S\pitchfork f$
to indicate that the map ${\cal E}(S,f):{\cal E}(S,A)\to {\cal E}(S,B)$
is surjective. If $\cal E$ has a terminal object $\top$,
the condition $u\pitchfork S$ is equivalent to the
condition $u\pitchfork t_S$, where $t_S$ is the map $S\to \top$.
If $\cal E$ has an initial object $\bot$, 
the condition $S\pitchfork f$ is equivalent to the
condition $i_S \pitchfork f$, where $i_S$ is the map $\bot \to S$.

\medskip

For any class of maps ${\cal M}\subseteq {\cal E}$,
we denote by ${}^\pitchfork\! {\cal M}$ (resp. ${\cal M}^\pitchfork$)
the class of maps having the left (resp. right) lifting property
with respect 
to every map in ${\cal M}$.
If ${\cal A}$ and ${\cal B}$
are two classes of maps in ${\cal E}$, we write ${\cal A}\pitchfork  {\cal B}$
to indicate that we have $a\pitchfork  b$ for every $a\in {\cal A}$
and $b\in {\cal B}$. 
Then
$$ {\cal A}\subseteq {}^\pitchfork {\cal B} \quad 
\Longleftrightarrow  \quad  {\cal A}\pitchfork  {\cal B}
\quad 
\Longleftrightarrow  \quad  {\cal B}\subseteq {\cal A}^\pitchfork .$$

\medskip

If $F:{\cal U}\leftrightarrow {\cal V}:G$
is a pair of adjoint functors, then for an arrow $f\in {\cal U}$
and an arrow $g\in  {\cal V}$ we have
$$f \pitchfork G(g) \quad 
\Longleftrightarrow  \quad  F(f)\pitchfork  g.$$

\medskip

\begin{definition} \label{defweakfact} We shall say that a pair $({\cal A},{\cal B})$ of classes 
of maps in a category ${\cal E}$ is a {\it weak factorisation system} 
if the following conditions are satisfied:
\begin{itemize}
\item {} every map $f\in {\cal E}$ admits a factorisation $f=pi$ 
with $i\in {\cal A}$ and $p\in {\cal B}$;
\item {} ${\cal B}={\cal A}^\pitchfork$ and ${\cal A}={}^\pitchfork {\cal B}$.
\end{itemize}
We call ${\cal A}$ is the {\it left class} and 
${\cal B}$ the {\it right class} of the weak factorisation system. 
\end{definition}

\begin{definition} \label{trivfib} We shall say that a map in a topos  is a {\it trivial fibration}
if it has the right lifting property
with respect to every monomorphism.\end{definition}

This terminology is non-standard but it is useful.  The trivial fibrations 
often coincide with the acyclic fibrations (which can be defined in any model category).

\begin{proposition} \label{trivfibfactsys} \cite{J2}
If ${\cal A}$ is the class of monomorphisms in  a 
topos and ${\cal B}$ is the class of trivial fibrations,
then the pair $({\cal A},{\cal B})$ is a weak factorisation system.
\end{proposition}

Recall that a map $u:A\to B$ in a category $\cal E$ is said to be
a {\it retract} of a map $f:X\to Y$ if $u$ is a retract of $f$
as objects of the category of arrows ${\cal E}^I$.
Recall that a map $u:A\to B$ is called
a {\it domain retract} of a map $v:C\to B$,
if $u$ is a retract of $v$ as objects of the category ${\cal E}/B$.
There is a dual notion of codomain retract.
The two classes of a weak factorisation system are closed under retracts.

\medskip

\begin{definition} \label{defsaturated} We shall say  that a class $\cal A$ of maps in a cocomplete
category $\cal E$
is {\it saturated} if it contains the isomorphisms
and is closed under composition, transfinite composition,
cobase change and codomain retracts. 
\end{definition}

The class ${}^\pitchfork{\cal M}$ is saturated
for any class ${\cal M}\subseteq \cal E$.
In particular, the class ${\cal A}$ of a weak factorisation system $({\cal A},{\cal B})$  in  $\cal E$
is saturated. Every class of maps ${\cal M}\subseteq  {\cal E}$
is contained in a smallest saturated class
${\overline {\cal M}}\subseteq  {\cal E}$ 
called the {\it saturated class generated} by ${\cal M}.$

\medskip

The following proposition is a special case of a more general result, see \cite{J2}:

\begin{proposition} \label{factorisationsystemgenerated} If $\Sigma$ 
is a set of maps in a presheaf category,
then the pair $({\overline \Sigma}, \Sigma^\pitchfork)$
is a weak factorisation system.
\end{proposition}

\medskip

We shall say that a functor of two variables
$$\odot :{\cal E}_1\times {\cal E}_2 \to {\cal E}_3$$
is  {\it divisible on the left} if the
functor $A\odot(-): {\cal E}_2 \to  {\cal E}_3$
admits a right adjoint $A\backslash (-): {\cal E}_3 \to  {\cal E}_2$
for every object $A\in {\cal E}_1$.
In this case we obtain a functor of two variables $(A,X)\mapsto A\backslash X$,
$${\cal E}_1^o\times {\cal E}_3 \to {\cal E}_2,$$
called the {\it left division functor}.
Dually, we shall say that $\odot$
is  {\it divisible on the right} if the
functor $(-)\odot B: {\cal E}_1 \to  {\cal E}_3$
admits a right adjoint $(-)/B: {\cal E}_3 \to  {\cal E}_1$
for every object $B\in {\cal E}_2$. 
In this case we obtain a functor of two variables
$(X,B)\mapsto X/B$,
$${\cal E}_3\times {\cal E}_2^o \to {\cal E}_1,$$
called the {\it right division functor}.
When the functor $\odot$ is divisible on both sides, 
there is a bijection between
the following three kinds of maps
$$A \odot B\to X, \quad   \quad B\to A \backslash X,\quad    \quad A\to X\slash B.$$
Hence the contravariant functors $A\mapsto A \backslash X$
and $B\mapsto B \backslash X$ are mutually right adjoint.
It follows that we have 
$$u\pitchfork (X/v)\quad \Leftrightarrow \quad v\pitchfork (u\backslash X).$$
for every map $u\in {\cal E}_1$
and every map $v\in  {\cal E}_2$.

\medskip

\noindent {\bf Remark}: If a functor of two variables $\odot :{\cal E}_1\times {\cal E}_2 \to {\cal E}_3$
is divisible on both sides, then so are
the left division functor ${\cal E}_1^o\times {\cal E}_3 \to {\cal E}_2$
and the right division functor  ${\cal E}_3\times {\cal E}_2^o \to {\cal E}_1.$
This is called a tensor-hom-cotensor situation in \cite{G}.

\medskip

Recall that a monoidal category ${\cal E}=( {\cal E},\otimes )$
is said to be {\it  closed} if the 
tensor product $\otimes $ is divisible on each side.
Let ${\cal E}=({\cal E},\otimes,\sigma)$ be a {\it symmetric} monoidal closed category,
with symmetry $\sigma:A\otimes B\simeq B\otimes A$.
Then the objects $X/A$ and $A\backslash X$ are canonicaly isomorphic;
we can identify them by adopting a common notation, for example $[A,X]$. 

\medskip

Recall that a category with finite products ${\cal E}$
is said to be  {\it cartesian closed} 
if the functor $A\times -: {\cal E} \to {\cal E}$ admits a right adjoint $(-)^A$
for every object $A\in {\cal E}$.
A cartesian closed category ${\cal E}$ is symmetric monoidal closed. 
Every presheaf category and more generally every
topos is cartesian closed.

\medskip

Let $\odot :{\cal E}_1\times {\cal E}_2 \to {\cal E}_3$
be a functor of two variables with values in a finitely cocomplete 
category ${\cal E}_3$.  If $u:A\to B$ is map in ${\cal E}_1$
and $v:S\to T$ is a map in ${\cal E}_2$,
we shall denote by $u   \odot' v$
the map 
$$A \odot  T  \sqcup_{A  \odot   S}  B \odot  S  \longrightarrow B  \odot    T
$$
obtained from the commutative square 
$$
\xymatrix{
A \odot S \ar[d] \ar[r]&  B  \odot  S  \ar[d]\\
A  \odot  T \ar[r]    &B  \odot    T .   }
$$
This defines a functor of two variables
$$\odot' :{\cal E}_1^I\times {\cal E}_2^I\to {\cal E}_3^I,$$
where ${\cal E}^I$ denotes the category of arrows
of a category ${\cal E}$. 

\medskip

In a topos, if $u:A\subseteq B$ and $v:S\subseteq T$ are inclusions of sub-objects
then the map $u\times' v$ 
is the inclusion of sub-objects 
$$(A\times T)\cup  (B  \times  S)\subseteq B\times T.$$

\medskip

Suppose now that a functor $\odot :{\cal E}_1\times {\cal E}_2 \to {\cal E}_3$
is divisible on both sides, that ${\cal E}_1$ and ${\cal E}_2$ are finitely complete 
and that  ${\cal E}_3$ is finitely cocomplete.
Then the functor $\odot' :{\cal E}_1^I\times {\cal E}_2^I\to {\cal E}_3^I$
is divisible on both sides.
If $u:A\to B$ is map in ${\cal E}_1$
and $f:X\to Y$ is a map in ${\cal E}_3$,
let us denote by $\langle u\backslash\ f\rangle$
the map 
$$B\backslash X \to B\backslash Y  \times_{ A\backslash Y }  A\backslash X$$
obtained from the
commutative square
$$\xymatrix{
B\backslash X \ar[r]\ar[d]&  A\backslash X \ar[d] \\
B\backslash Y   \ar[r]&  A\backslash Y .
 }
 $$
Then the functor  $f\mapsto \langle u\backslash f\rangle$
is right adjoint to the functor $v\mapsto u\odot' v$.
Dually, if $v:S\to T$ is map in ${\cal E}_2$
and $f:X\to Y$ is a map in ${\cal E}_3$,
we shall denote by $\langle f/v\rangle$
the map 
$$X/T \to Y/T \times_{ Y/S}  X/S$$
obtained from the
commutative square
$$\xymatrix{
X/T \ar[r]\ar[d]&  X/S\ar[d] \\
Y/T  \ar[r]&  Y/S .
 }
$$
The functor  $f\mapsto \langle f\backslash v\rangle$
is right adjoint to the functor $u\mapsto u\odot' v$.

\medskip

The verification of the following result is left to the reader.  See \cite{J2}.

\begin{proposition}  \label{lift2A} 
Let $\odot :{\cal E}_1\times {\cal E}_2 \to {\cal E}_3$
be a functor of two variables divisible on both sides,
where ${\cal E}_i$ is a finitely bicomplete category for $i=1,2,3$.
If  $u\in {\cal E}_1$, $v\in {\cal E}_2$ and $f\in {\cal E}_3$, then  
$$(u  \odot' v) \pitchfork f \quad \Longleftrightarrow \quad u\pitchfork \langle f/v \rangle \quad 
\Longleftrightarrow  \quad v\pitchfork \langle u\backslash f \rangle.$$
\end{proposition}

\medskip

Let ${\cal E}=({\cal E},\otimes)$ be a bicomplete symmetric monoidal closed category.
Then the operation $\otimes'$ gives the category ${\cal E}^I$
 the structure of a symmetric monoidal closed category. 
 If $u$ and $f$
are two maps in ${\cal E}$, then
the maps $\langle f/u\rangle$ and $\langle u\backslash f\rangle $
are canonically isomorphic. We shall
identify them by adopting a common notation $\langle u, f\rangle$.
 If  $u:A\to B$,  $v:S\to T$ 
and  $f:X\to Y$ are three maps in $\cal E$, then  
$$(u  \otimes' v) \pitchfork f \quad \Longleftrightarrow \quad u\pitchfork \langle v,f \rangle \quad 
\Longleftrightarrow  \quad v\pitchfork \langle u, f \rangle.$$

\bigskip

We now recall the notion of a Quillen model category:

\begin{definition} \cite{Q} \label{defmodelcat} Let ${\cal E}$ be a finitely bicomplete category.
A {\it model structure} on ${\cal E}$
is a triple $({\cal C},{\cal W},{\cal F})$
of classes of maps in ${\cal E}$ 
satisfying   
the following conditions:
\begin{itemize}
\item {} {\it (``three-for-two")} if two of the three maps 
$u:A\to B$, $v:B\to C$ and $vu:A\to C$ belong to ${\cal W}$, then so does the third; 
\item {} the pair $( {\cal C}\cap {\cal W},{\cal F})$ is a weak factorisation
system; 
\item {} the pair $({\cal C},{\cal F}\cap {\cal W})$ is
a weak factorisation system. 
\end{itemize}
\end{definition}

These conditions imply that ${\cal W}$
is closed under retracts by \ref{Wclosedunderretract} below.
A category ${\cal E}$ equipped with
a model structure is called a  {\it model category}.
A map in ${\cal C}$ is called a {\it cofibration}, a map in ${\cal F}$ a {\it fibration}
and a  map in ${\cal W}$ a  {\it weak equivalence}.
A map in ${\cal W}$ is also said to be  {\it acyclic}.
An object $X\in {\cal E}$ is {\it fibrant} if the map
$X\to 1$ is a fibration, where $1$ is the terminal object of 
${\cal E}$. Dually, an object $A\in {\cal E}$ is  {\it cofibrant} if the
map $0\to A$ is a cofibration, where $0$
is the initial object of ${\cal E}$.

\medskip

Any two of the three classes of a model structure $({\cal C},{\cal W},{\cal F})$
determine the third.

\medskip

A model structure is said to be {\it left proper} if the cobase change
of an acyclic map along a cofibration is acyclic.
Dually, a model structure is said to be {\it right proper} if the base change 
of an acyclic map
along a fibration is acyclic.
A model structure is {\it proper} if it is both left and right proper.

\begin{proposition} \label{Wclosedunderretract} \cite {JT2} The class ${\cal W}$ of a model structure is closed under retracts.\end{proposition}

\noindent
{\bf Proof}: Observe first that the class ${\cal F}\cap {\cal W}$ 
is closed under retracts since the pair $({\cal C},{\cal F}\cap {\cal W})$ 
is a weak factorisation system.
Suppose now that a map $f:A\to B$ is a retract of a map $g:X\to Y$ in ${\cal W}$.
Let us show that $f\in {\cal W}$.
By definition, we have a commutative diagram,
 $$
\xymatrix{
A \ar[d]_f \ar[r]^s & X\ar[r]^t \ar[d]^{g}& A\ar[d]^f \\
B \ar[r]_u & Y\ar[r]^v & B }
$$
where $gf=1_X$ and $vu=1_B$.
Let us first consider the case where $f$ is a fibration.
In this case, let us choose a factorisation $g=qj:X\to Z\to Y$
with $j\in{\cal C} \cap {\cal W}$ and $q\in {\cal F}$.
We have $q\in {\cal F} \cap {\cal W}$ by three-for-two, since $g\in {\cal W}$.
The square 
$$
\xymatrix{
X \ar[d]_j \ar[r]^t &  A\ar[d]^f \\
Z \ar[r]^{vq} &  B }
$$
has a diagonal filller $d:Z\to A$, since $f$ is a fibration.
We then have a commutative diagram,
$$
\xymatrix{
A \ar[d]_f \ar[r]^{js} & Z\ar[r]^d \ar[d]_{q}& A\ar[d]^f \\
B \ar[r]^u & Y\ar[r]^v & B.}
$$
Thus, $f$ is a retract of $q$, since $d(js)=ts=1_A$. 
This shows that $f\in  {\cal W}$ since  $q\in {\cal F} \cap {\cal W}$.
In the general case,
let us choose a factorisation $f=pi:A\to E\to B$
with $i\in{\cal C} \cap {\cal W}$ and $p\in {\cal F}$
By taking a pushout, we obtain a commutative diagram
$$
\xymatrix{
A \ar[d]_i \ar[r]^s & X\ar[r]^t \ar[d]^{i_2}& A\ar[d]^i \\
E \ar[d]_p \ar[r]^(0.4){i_1} & E\sqcup_A X\ar[r]^r \ar[d]^{k}& E\ar[d]^p \\
B \ar[r]^u & Y\ar[r]^v & B, }
$$
where $ki_2=g$ and $ri_1=1_E$. The map $i_2$ is a cobase
change of the map $i$. Thus, $i_2\in {\cal C} \cap {\cal W}$
since $i\in  {\cal C} \cap {\cal W}$. Thus, $k\in {\cal W}$
by three-for-two since $g=ki_2\in {\cal W}$ by hypothesis.
Thus, $p\in {\cal W}$ by the first part since $p\in  {\cal F}$.
Thus $f=pi\in  {\cal W}$ since $i\in  {\cal W}$.
\QED

\medskip

The {\it homotopy category} of a model category ${\cal E}$
is defined to be the category of fractions $Ho({\cal E})={\cal W}^{-1}{\cal E}$.
We shall denote by $[u]$ the image of an arrow $u\in {\cal E}$
by the canonical functor ${\cal E}\to Ho({\cal E})$.
The arrows $[u]$ is invertible iff $u$ is a weak equivalence
by a result in \cite{Q}.

\medskip

Let ${\cal E}_f$ (resp. ${\cal E}_c$) be  the full sub-category of fibrant (resp. cofibrant)
objects of ${\cal E}$ and let us put ${\cal E}_{fc}= {\cal E}_f\cap {\cal E}_c$.
Let us put $Ho({\cal E}_{f})= {\cal W}_{f}^{-1}{\cal E}_{f}$ where ${\cal W}_f={\cal W}\cap {\cal E}_f$
and similarly for $Ho({\cal E}_{c})$ and $Ho({\cal E}_{fc})$.
.
Then the diagram of inclusions 
$$
\xymatrix{
{\cal E}_{fc}  \ar[d] \ar[r] & {\cal E}_f \ar[d] \\
{\cal E}_{c}  \ar[r] & {\cal E} }
$$
induces a diagram of equivalences of categories 
$$
\xymatrix{
Ho({\cal E}_{fc}) \ar[d] \ar[r] & Ho({\cal E}_f) \ar[d]\\
Ho({\cal E}_{c})  \ar[r] & Ho({\cal E}).}
$$

\medskip

A {\it fibrant replacement} of an object $X\in {\cal E}$ 
is a weak equivalence $X\to RX$ with codomain a fibrant object.
Dually, a {\it cofibrant replacement} of $X$ is a
 weak equivalence $LX\to X$ with domain a cofibrant object.

\bigskip

Recall from \cite{Ho} that a cocontinuous functor $F:{\cal U}\to {\cal V}$
between two model categories
is said to be a {\it left Quillen functor} if  it takes
a cofibration to a cofibration and an acyclic
cofibration to and acyclic cofibration. 
Dually, a continuous functor $G:{\cal V} \to {\cal U}$ between two model categories
is said to be a {\it right Quillen functor} if it takes
a fibration to a fibration and an acyclic
fibration to an acyclic fibration.

\begin{proposition}   \cite {Q} \label{Quillenpair1} Let $F:{\cal U}\leftrightarrow {\cal V}:G$
be an adjoint pair of functors between two model categories.
Then $F$ is a left Quillen functor iff  $G$ is a right Quillen functor.
\end{proposition}

The adjoint pair $(F,G)$ is said to be a {\it Quillen pair}
if the conditions of \ref{Quillenpair1} are satisfied.

\medskip

The following lemma is due to Ken Brown, see \cite {Ho} and \cite{JT2}.

\begin{lemma}  \label{KBlemma0} 
Let $\cal E$ be a model category and
$F:{\cal E}\to {\cal D}$ be a functor with values 
in a category equipped with a class of a weak equivalences
${\cal W}'$ which satisfies three-for-two.
If $F$ takes an acyclic cofibration between
cofibrant objects to a weak equivalence,
then it takes a weak equivalence
between cofibrant objects to a weak equivalence.
\end{lemma}

\begin{corollary}  \label{KBlemma} 
A left Quillen functor
takes a weak equivalence between cofibrant objects
to a weak equivalence. 
\end{corollary}

The following result is due to Reedy \cite{Ree}.
See \cite {Hi}  and \cite{JT2}.

\begin{proposition}  \label{Leftproper0} The cobase change along
a cofibration of a weak
equivalence between cofibrant objects is a weak equivalence.
\end{proposition}

\begin{corollary} \label{Leftproper} If every object of model category
is cofibrant
then the model structure is left proper.
\end{corollary}

\begin{lemma} \label{lift3}  In a model category, a cofibration is acyclic iff 
it has the left lifting property with respect
to  every fibration between fibrant objects.
\end{lemma}

\noindent {\bf Proof}: The necessity is clear. 
Conversely, let us suppose that a cofibration
$u:A\to B$ has the left lifting property with respect
to  every fibration between fibrant objects.  We shall prove that
$u$ is acyclic.
For this, let us choose a fibrant replacement $j:B\to B'$ 
of the object $B$
together with a factorisation of the composite $ju:A\to B'$
as a weak equivalence $i:A\to A'$ followed by a fibration $p:A'\to B$.
The square
$$
\xymatrix{
A  \ar[d]_u \ar[r]^{i}& A' \ar[d]^p\\
B \ar[r]^{j}   & B'     }
$$
has a diagonal filler $d:B\to A'$ 
since $p$ is a fibration between fibrant objects.
The arrows $i$ and $j$ are invertible in the homotopy
category since they are acyclic.
The relations $pd=j$ and  $du=i$
then implies that $d$ is invertible in the homotopy
category. It thus acyclic \cite{Q}.
It follows by three-for-two that $u$ is acyclic.
\QED

\begin{proposition} \label{Quillenpair2} An adjoint pair of functors
$F:{\cal U}\leftrightarrow {\cal V}:G$
between two model categories
is a Quillen pair iff
 the following two conditions are satisfied:
 \begin{itemize}
\item $F$ takes a cofibration to a cofibration;
\item $G$ takes a fibration between fibrant objects
to a fibration. 
\end{itemize}
\end{proposition}

\noindent{\bf Proof}: The necessity is obvious.
Let us prove the sufficiency.
For this it suffices to show that $F$
is a left Quillen functor by \ref{Quillenpair1}.
Thus we show  that $F$ takes an acyclic
cofibration $u:A\to B$ to an acyclic  cofibration $F(u):F(A)\to F(B)$.
But $F(u)$ is acyclic iff it has the left
lifting property with respect to every fibration between
fibrant objects $f:X\to Y$ by Lemma \ref{lift3}.
But the condition $F(u)\pitchfork f$
is equivalent to the condition
$u\pitchfork G(f)$ by the  adjointness  $F\dashv G$.
We have $u\pitchfork G(f)$
since $ G(f)$
is a fibration by (ii).
This proves that we have $F(u)\pitchfork f$.
Thus, $F(u)$ is acyclic. 
\QED

A left Quillen functor $F:{\cal U}\to {\cal V}$
induces a functor $F_c:{\cal U}_c\to {\cal V}_c$
hence also a functor $Ho(F_c):Ho({\cal U}_c) \to  Ho({\cal V}_c)$
by Proposition \ref{KBlemma}.
A {\it left derived functor} is a functor 
$$F^L:Ho({\cal U})\to Ho({\cal V})$$
for which 
the following diagram of functors commutes up to 
isomorphism,
$$
\xymatrix{
Ho({\cal U}_c)  \ar[d] \ar[r]^{Ho(F_c)} & Ho({\cal V}_c)\ar[d]\\
Ho({\cal U})  \ar[r]^{F^L} & Ho({\cal V}), }
$$
The functor $F^L$ is unique up to a canonical isomorphism.
It can be computed as follows.
For each object $A\in {\cal U}$, we can choose a cofibrant replacement 
$\lambda_A:LA\to A$, with $\lambda_A$ an acyclic fibration.
We can then choose for each arrow $u:A\to B$ an arrow
$L(u):LA\to LB$ such that $u\lambda_A=\lambda_B  L(u)$,
$$
\xymatrix{
LA \ar[d]_{L(u)} \ar[r]^{\lambda_A}& A\ar[d]^u\\
LB \ar[r]^{\lambda_B}   & B . }
$$
Then 
$$F^L([u])=[F(L(u))] :FLA\to FLB.$$

\bigskip

A right Quillen functor $G:{\cal V} \to {\cal U}$ 
induces a functor $G_f:{\cal V}_f\to {\cal U}_f$
hence also a functor $Ho(G_f):Ho({\cal V}_f) \to  Ho({\cal U}_f)$
by Proposition \ref{KBlemma}.
The {\it right derived functor} is a functor 
$$G^R:Ho({\cal V})\to Ho({\cal U})$$
for which 
the following diagram of functors commutes up to a canonical
isomorphism,
$$
\xymatrix{
Ho({\cal V}_f)  \ar[d] \ar[r]^{Ho(G_f)} & Ho({\cal U}_f)\ar[d]\\
Ho({\cal V})  \ar[r]^{G^R} & Ho({\cal U}). }
$$
The functor $G^R$ is unique up to a canonical isomorphism.
It can be computed as follows.
For each object $X\in {\cal V}$ let us choose a fibrant replacement
$\rho_X:X\to RX$, with $\rho_X$ an acyclic cofibration.
We can then choose for each arrow $u:X\to Y$ an arrow
$R(u):RX\to RY$ such that $R(u)\rho_X=\rho_Y u$,
$$
\xymatrix{
X \ar[d]_{u} \ar[r]^{\rho_X}& RX\ar[d]^{R(u)}\\
Y \ar[r]^{\rho_Y}   & RY . }
$$
Then 
$$G^R([u])=[G(R(u))]:GRX\to GRY.$$

\bigskip

A Quillen pair of adjoint functors $F:{\cal U}\leftrightarrow {\cal V}:G$
induces
a pair of adjoint functors
$$F^L:Ho({\cal U})\leftrightarrow Ho({\cal V}):G^R.$$
If $A\in {\cal U}$ is cofibrant, the adjunction unit $A\to G^RF^L(A)$
is obtained by composing the maps $A\to GFA\to GRFA$,
where $FA\to  RFA$ is a fibrant replacement of $FA$.
If $X\in {\cal V}$ is fibrant, the adjunction counit $F^LG^R(X)\to X$
is obtained by composing the maps $FLGX \to FGX\to X$,
where $LGX\to  GX$ is a cofibrant replacement of $GX$.

\begin{definition} \label{Localisation1} We shall say that
a Quillen pair
$F:{\cal U}\leftrightarrow {\cal V}:G$
is a {\it homotopy localisation} ${\cal U}\to {\cal V}$  if the right derived functor $G^R$
is full and faithful.
Dually, we shall say that the pair
$(F,G)$ is a {\it homotopy colocalisation} ${\cal V}\to {\cal U}$
if the left derived functor $F^L$
is full and faithful.
\end{definition}

\begin{proposition} \label{Localisation2} The following conditions on a Quillen pair $F:{\cal U}\leftrightarrow {\cal V}:G$ are equivalent:
\begin{itemize}
\item The pair $(F,G)$ is a homotopy localisation ${\cal U}\to {\cal V}$;
\item The map $FLGX\to X$
is a weak equivalence for every fibrant object $X\in   {\cal V}$, where $LGX\to GX$
denotes a cofibrant replacement of $GX$;
\item The map $FLGX\to X$
is a weak equivalence for every fibrant-cofibrant object $X\in   {\cal V}$, where $LGX\to GX$
denotes a cofibrant replacement of $GX$.
\end{itemize}
\end{proposition}

\noindent {\bf Proof}: 
The functor $G^R$ is full and faithful iff the counit of the adjunction $F^L\dashv G^R$
is an isomorphism. But if $X\in {\cal V}$ is fibrant, this counit is obtained by composing the maps $FLGX \to FGX\to X$, where $LGX\to  GX$ is a cofibrant replacement of $GX$.
This proves the equivalence (i)$\Leftrightarrow$(ii). 
The implication (ii)$\Rightarrow$(iii) is obvious.
Let us prove the implication (iii)$\Rightarrow$(ii). 
For every fibrant objet $X$, there is a an acyclic fibration $p:Y\to X$
with domain a cofibrant object $Y$. 
The map $Gp:GY\to GX$
is an acyclic fibration, since $G$ is a right Quillen functor.
Let $q:LGY\to GY$ be a cofibrant replacement of $GY$.
Then  the map $FLGY\to FGY \to Y$
is a weak equivalence by assumption, since $Y$ is fibrant-cofibrant.
But the composite $G(p)q: LGY\to GY\to GX$ is a cofibrant
replacement of $GX$, since $G(p)$ is a weak equivalence.
Moreover, the composite $FLGY\to FGX \to X$ is a weak equivalence,
since $p$ is a weak equivalence and
the following diagram commutes 
$$
\xymatrix{
FLGY \ar[dr]  \ar[r] & FGY \ar[d] \ar[r]  & Y \ar[d]^{p} \\
 & FGX \ar[r] & X .
 }
$$
This proves that condition (ii) is satisfied for a cofibrant
replacement of $GX$.
\QED

\begin{proposition} \label{Localisation3}  If $F:{\cal U}\leftrightarrow {\cal V}:G$
is a homotopy localisation, then the right adjoint $G$ preserves and reflects weak equivalences between fibrant objects.
\end{proposition}

\noindent {\bf Proof}: The functor $G^R$
is equivalent to the functor $Ho(G_f):Ho({\cal V}_f) \to  Ho({\cal U}_f)$
induced by the functor $G$. Thus, $Ho(G_f)$
is full and faithful since $G^R$ is full and faithful.
This proves the result
since a full and faithful functor is conservative.
\QED

\begin{proposition} \label{Localobject2A}  Let $M_i=({\cal C}_i,{\cal W}_i, {\cal F}_i)$ ($i=1,2$)
be two model structures on a category ${\cal E}$.
Suppose that ${\cal C}_1\subseteq {\cal C}_2$
and ${\cal W}_1\subseteq {\cal W}_2$.
Then the identity
functor ${\cal E} \to {\cal E}$ is a homotopy localisation $M_1\to M_2$.
\end{proposition}

\begin{definition} \label{DefBousfieldLocalisationA}  
Let $M_i=({\cal C}_i,{\cal W}_i, {\cal F}_i)$ ($i=1,2$)
be two model structures on a category ${\cal E}$.
If ${\cal C}_1= {\cal C}_2$
and ${\cal W}_1\subseteq {\cal W}_2$, we shall say that 
$M_2$
is a {\it Bousfield localisation} of $M_1$.
\end{definition}

\begin{proposition} \label{Localisation5A}  Let $M_2=({\cal C}_2,{\cal W}_2, {\cal F}_2)$ be
a Bousfield localisation of a model structure $M_1=({\cal C}_1,{\cal W}_1, {\cal F}_1)$
on a category ${\cal E}$. 
Then a map between $M_2$-fibrant objects
is a $M_2$-fibration 
iff it is a $M_1$-fibration.
\end{proposition}

\noindent {\bf Proof}: 
By hypothesis, we have ${\cal C}_1= {\cal C}_2$ and
${\cal W}_1\subseteq {\cal W}_2$. It follows that we have
${\cal F}_2\cap {\cal W}_2= {\cal F}_1\cap {\cal W}_1$ and $ {\cal F}_2\subseteq {\cal F}_1$.
Let $f:X\to Y$ be a map  between two $M_2$-fibrant objects.
Let us show that $f$ is a $M_2$-fibration 
iff it is a $M_1$-fibration.
The implication ($\Rightarrow$) is clear, since $ {\cal F}_2\subseteq {\cal F}_1$.
Conversely, if $f\in {\cal F}_1$, let us show that $f\in {\cal F}_2$.
Let us choose a factorisation $f=pi:X\to Z\to Y$
with $i\in {\cal C}_2\cap {\cal W}_2$
and $p\in {\cal F}_2$. We have $i\in {\cal W}_1$
by Proposition \ref{KBlemma}, since
the identity functor is a right Quillen
functor $M_2\to M_1$
and since $i$ is a map between $M_2$-fibrant objects.
Thus, $i\in {\cal W}_1\cap {\cal C}_1$, since ${\cal C}_1={\cal C}_2$.
Hence the square
$$
\xymatrix{
 X  \ar[d]_i \ar[r]^{id} & X \ar[d]^f \\
E \ar[r]_p     &Y      }
$$
has a diagonal filler,
making $f$ a retract of $p$ and therefore $f\in {\cal F}_2$.
\QED

A Quillen pair $(F,G)$ is said to be a {\it Quillen equivalence}
if the adjoint pair 
$(F^L,G^R)$ is an equivalence of categories.

\begin{proposition} \label{Quillenequiv} A Quillen pair $F:{\cal U}\leftrightarrow {\cal V}: G$
is a Quillen equivalence iff the following equivalent conditions are satisfied:
\begin{itemize}
\item The pair $(F,G)$ is a both a homotopy localisation and colocalisation;
\item The pair $(F,G)$ is a homotopy localisation and the functor $F$ reflects weak equivalences between 
cofibrant objects;
\item The pair $(F,G)$ is a homotopy colocalisation and the functor $G$ reflects weak equivalences between fibrant objects;
\end{itemize}
\end{proposition}

The composite of two adjoint pairs 
$$F_1:{\cal E}_1\leftrightarrow {\cal E}_2:G_1 \quad {\rm and} \quad
F_2:{\cal E}_2\leftrightarrow {\cal E}_3:G_2$$
is an adjoint pair
$F_2F_1:{\cal E}_1\leftrightarrow {\cal E}_3:G_1G_2.$


\begin{proposition} [Three-for-two, \cite{Ho}] \label{Quillenequiv3for2}  The composite
of two Quillen pairs 
$(F_1,G_1)$ and $(F_2,G_2)$ is a Quillen pair $(F_2F_1,G_1G_2)$.
Moreover, if two of the pairs $(F_1,G_1)$, $(F_2,G_2)$ and $(F_2F_1,G_1G_2)$
are Quillen equivalences, then so is the third.
\end{proposition}

\medskip

\begin{definition} \cite{Ho} \label{defQuillen2v} We shall say that
a functor of two variables between three model categories
 $$\odot :{\cal E}_1\times {\cal E}_2 \to {\cal E}_3$$ 
is a {\it left Quillen functor} if it is cocontinuous in each variable
and the following conditions are satisfied:
\begin{itemize}
\item {}  $u \odot' v $ is a cofibration 
if $u\in  {\cal E}_1$
and $v\in  {\cal E}_2$ are cofibrations; 
\item {}  $u \odot' v $ is an acyclic cofibration 
if $u\in  {\cal E}_1$
and $v\in  {\cal E}_2$ are cofibrations and if $u$ or $v$ is acyclic.
\end{itemize}
Dually, we shall say that $\odot$
is a {\it right Quillen functor}
if the opposite functor  $\odot^o :{\cal E}_1^o\times {\cal E}_2^o \to {\cal E}_3^o$
is a left Quillen functor.
\end{definition}

\medskip

\begin{proposition} \label{propQuillen2vto1v} Let $\odot:{\cal E}_1\times {\cal E}_2 \to {\cal E}_3$ be a left Quillen functor of two variables between three model 
categories.  If $A\in {\cal E}_1$ is cofibrant, then
the functor $B\mapsto A\odot B$ is a left Quillen functor ${\cal E}_2 \to  {\cal E}_3$.
\end{proposition}

\noindent{\bf Proof}: If $A\in {\cal E}_1$ is cofibrant, then the map $i_A:\bot \to A$ is a cofibration,
where $\bot$ is the initial object.
If $v:S\to T$ is a map in ${\cal E}_2$, then we have 
$A\odot v=i_A\odot' v$.
Thus, $A\odot v$ is a cofibration if $v$ is a cofibration
and $A\odot v$ is acyclic if moreover $v$ is acyclic. \QED

\begin{proposition}  Let $\odot :{\cal E}_1\times {\cal E}_2 \to {\cal E}_3$
be a functor of two variables between three model categories. 
If the functor $\odot$ is divisible on the left, then it is
a left Quillen functor iff the corresponding left division functor ${\cal E}_1^o\times {\cal E}_3 \to {\cal E}_2$ is a right Quillen functor.
Dually, if the functor $\odot$ is divisible on the right, then it is
a left Quillen functor iff the corresponding right division functor 
${\cal E}_3 \times {\cal E}^o_2 \to {\cal E}_1$ is a right Quillen functor.
\end{proposition}

\medskip

\begin{proposition} \label{propQuillen2v} 
Let $\odot:{\cal E}_1\times {\cal E}_2 \to {\cal E}_3$ be a functor of two variables,
cocontinuous in each, 
between three model categories. Suppose that 
the following three conditions are satisfied:
\begin{itemize}
\item{} If $u\in {\cal E}_1$ and $v\in {\cal E}_2$ are cofibrations, then so is $u \odot' v$;
\item{} the functor $(-)\odot B$ preserves acyclic cofibrations for every object $B\in  {\cal E}_2$;
\item{} the functor $A\odot (-)$ preserves acyclic cofibrations for every object $A\in  {\cal E}_1$.
\end{itemize}
Then $\odot $ is a left Quillen functor.
\end{proposition}

\noindent{\bf Proof}: Let $u:A\to B$ be a cofibration in  ${\cal E}_1$
and $v:S\to T$ be a cofibration in  ${\cal E}_2$.
Let us show that $u \odot' v$ is acyclic if $u$ or $v$ is acyclic.
We only consider the case where $v$ is acyclic.
Consider the commutative diagram 
$$
\xymatrix{
A\odot S \ar[d]_{A\odot v} \ar[r]^{u\odot S}&  B  \odot  S  \ar[d]_{i_2}  \ar[dr]^{B\odot v} & \\
A  \odot  T \ar[r]^{i_1}    & Z  \ar[r]^(0.4){u \odot' v}    & B  \odot   T  \\
}
$$
where $Z=A  \odot  T\sqcup_{A \odot S} B  \odot  S$
and where  $(u \odot' v)i_1= u\odot T$. The map $A\odot v$ is
an acyclic cofibration since $v$ is an acyclic cofibration.
Similarly for the map $B\odot v$.
It follows that $i_2$ is an acyclic cofibration by cobase change.
 Thus, $u \odot' v$ is acyclic by three-for-two since 
 $(u \odot' v)i_2=B\odot v$  is acyclic. 

\medskip

\begin{definition} \cite {Ho} A 
model structure $({\cal C},{\cal W},{\cal F})$ on a monoidal closed category 
 ${\cal E}=({\cal E},\otimes)$ is said to be {\it monoidal} if the tensor product 
$\otimes:{\cal E} \times {\cal E} \to {\cal E}$
 is a left Quillen functor of two variables
 and if the unit object of the tensor product is
cofibrant. 
\end{definition}

\medskip

In a monoidal closed model category, if $f$ is a fibration 
then so are the maps $\langle u\backslash f\rangle $ and $\langle f/u \rangle$
for any cofibration $u$. Moreover, the
fibrations $\langle u\backslash f\rangle $ and $\langle f/u \rangle$
are acyclic if the cofibration $u$ is acyclic or the fibration $f$ is acyclic. 

\medskip

\begin{definition}\label{defcartclosedmod} We shall say that a
model structure $({\cal C},{\cal W},{\cal F})$ on a cartesian closed category ${\cal E}$
 is {\it cartesian closed} if the cartesian product 
 $\times:{\cal E} \times {\cal E} \to {\cal E}$
 is a left Quillen functor of two variables
 and if the terminal object  1 is cofibrant.
 \end{definition}

In a cartesian closed model category,  if 
$f$ is a fibration and $u$ is a cofibration,
then the map $\langle u, f\rangle $
is a fibration, which is acyclic if $u$ or $f$ is acyclic.

\medskip

We recall a few notions of enriched category theory \cite{K}.
Let ${\cal V}=({\cal V},\otimes,\sigma)$ a  bicomplete symmetric monoidal closed category. 
A category enriched over $\cal V$ is called  a $\cal V$-{\it category}.
If ${\cal A}$ and ${\cal B}$ are $\cal V$-categories,
there is a notion of a {\it strong} functor $F:{\cal A}\to {\cal B}$;
it is an ordinary functor equipped with a {\it strength} which
is a natural transformation ${\cal A}(X,Y)\to {\cal B}(FX,FY)$ preserving
composition and units. A natural transformation $\alpha:F\to G$
between strong functors $F,G:{\cal A}\to {\cal B}$ is said to be {\it strong}
if the following square commutes
$$
\xymatrix{
{\cal A}(X ,Y) \ar[d] \ar[rr] && {\cal B}(GX,GY)  \ar[d]^{ {\cal B}(\alpha_X,GY)}  \\
{\cal B}(FX,FY) \ar[rr]^{  {\cal B}(FX,\alpha_Y)}    &&{\cal B}(FX, GY) . } 
$$
for every pair of objects $X,Y\in {\cal A}$. 
A {\it strong adjunction} $\theta:F\dashv G$
between strong functors $F:{\cal A}\to {\cal B}$ and $G:{\cal B}\to {\cal A}$
is a strong natural isomorphism 
$$\theta_{XY}:{\cal A}(FX ,Y)\to {\cal B}(X ,GY).$$
A strong functor $G:{\cal B}\to {\cal A}$ has
a strong left adjoint iff it has an ordinary left adjoint $F:{\cal A}\to {\cal B}$
and the
map 
$$
\xymatrix{  {\cal B}(FX ,Y) \ar[r]  & {\cal A}(GFX,GY)  \ar[rr]^{{\cal A}(\eta_X,GY)} &&  {\cal A}(X,GY)  }
$$
obtained by composing with the unit $\eta_X$ of the adjunction is an isomorphism for every pair of
 objects
$X\in  {\cal A}$ and $Y\in  {\cal B}$.
Recall that a $\cal V$-category ${\cal E}$ is said to admit {\it tensor products} 
if the functor $Y\mapsto {\cal E}(X,Y)$ admits
a strong left adjoint $A\mapsto A\otimes X$ for every object $X\in {\cal E}$.
A  $\cal V$-category  is said to be  (strongly) {\it cocomplete}
if it cocomplete as an ordinary category and 
if it admits tensor products.
These notions can be dualised. 
A $\cal V$-category ${\cal E}$ is said to admit {\it cotensor products} 
if the opposite $\cal V$-category ${\cal E}^o$
admits tensor products.  
This means that the (contravariant) functor $X\mapsto {\cal E}(X,Y)$ admits
a strong right adjoint $A\mapsto Y^{[A]}$ for every object $Y\in {\cal E}$.
A $\cal V$-category is said to be (strongly) {\it complete} 
if it complete as an ordinary category and if it
admits cotensor products.
We shall say that a $\cal V$-category  is (strongly) {\it bicomplete}
if it is both $\cal V$-complete and cocomplete.

\begin{definition} \cite{Q} Let ${\cal E}$ be a strongly bicomplete simplicial category.
We shall say that a
model structure $({\cal C},{\cal W},{\cal F})$ on ${\cal E}$ is {\it simplicial} if 
the tensor product $$\otimes:{\bf S}\times {\cal E}\to {\cal E}$$
is a left Quillen functor of two variables, where ${\bf S}$ is equipped
with the model structure $({\cal C}_0,{\cal W}_0,{\cal F}_0)$ of \ref{classicalmodel}.
\end{definition}

A simplicial category equipped with a simplicial model structure
is called a {\it simplicial model category}.

\medskip

\begin{proposition}  \label{simpmodcatequiv} Let $\cal E$
be a simplicial model category.
Then a map between cofibrant objects
$u:A\to B$ is acyclic iff
the map of simplicial sets
$${\cal E}(u,X):{\cal E}(B,X)\to {\cal E}(A,X)$$
is a weak homotopy equivalence
for every fibrant object $X$.
\end{proposition}

\noindent{\bf Proof}: The functor $A\mapsto  {\cal E}(A,X)$
takes
an (acyclic) cofibration to an (acyclic) Kan fibration
if $X$ is fibrant.
It then follows by Proposition \ref{KBlemma}
that it takes 
an acyclic map between cofibrant objects to an acyclic map.
Conversely, let $u:A\to B$ be a map between cofibrant objects in ${\cal E}$.
 If the map ${\cal E}(u,X):{\cal E}(B,X)\to {\cal E}(A,X)$
 is a weak homotopy equivalence
for every fibrant object $X$, 
let us show that $u$ is acyclic.
Let us first suppose that  $A$ and $B$
are fibrant. Let ${\cal E}_{cf}$
be the full subcategory of fibrant and cofibrant objects of ${\cal E}$.
We shall prove
that $u$ is acyclic by showing that
$u$ is invertible in the homotopy category $Ho({\cal E}_{cf})$.
But if $S,X\in {\cal E}_{cf}$, then we have
$Ho({\cal E}_{cf})(S,X)=\pi_0{\cal E}(S,X)$ by \cite{Q}.
Hence the map $Ho({\cal E}_{cf})(u,X): Ho({\cal E}_{cf})(B,X)\to Ho({\cal E}_{cf})(A,X)$
is equal to the map $\pi_0{\cal E}(u,X):\pi_0{\cal E}(B,X)\to \pi_0{\cal E}(A,X)$.
But the map $\pi_0{\cal E}(u,X)$ is bijective since 
the map ${\cal E}(u,X)$
is a weak homotopy equivalence.
This shows that the map $Ho({\cal E}_{cf})(u,X)$
is bijective for every $X\in {\cal E}_{cf}$.
It follows by the Yoneda lemma
that $u$ is invertible in $Ho({\cal E}_{cf})$.
Thus, $u$ is acyclic by \cite{Q}. 
In the general case,
let us choose a fibrant replacement  $i_A:A\to A'$ with
$i_A$ an acyclic cofibration.
Similarly, let us choose a fibrant replacement  $i_B:B\to B'$ with
$i_B$ an acyclic cofibration.
Then there exists a map  $u':A'\to B'$ 
such that $u'i_A=i_Bu$. 
We then have a commutative square 
$$
\xymatrix{
{\cal E}(B',X) \ar[d] \ar[r]&  {\cal E}(A',X)  \ar[d]\\
{\cal E}(B,X) \ar[r]    &{\cal E}(A,X)   }
$$
for every object $X$.
If $X$ is fibrant, then the vertical maps of the square  
are weak homotopy equivalences 
by the first part of the proof.
Hence also the map ${\cal E}(u',X):{\cal E}(B',X)\to {\cal E}(A',X)$
by  three-for-two.
This shows that $u'$ is an acyclic map since $A'$ and $B'$ are fibrant.
 It follows by three-for-two that $u$ is acyclic. 
\QED

\bigskip

Let $\cal E$ be a bicomplete category.
The  {\em box product} of a simplicial set $A$ by an object  $B\in {\cal E}$
is defined to be the simplicial object $A \Box B\in [\Delta^{o}, {\cal E}]$
obtained by putting
$$(A \Box B)_n =A_n\times B$$
for every $n\geq 0$, 
where $A_n\times B$ denotes the coproduct of $A_n$ copies of the object $B$.
The  functor 
$ \Box :{\bf S} \times {\cal E} \to [\Delta^{o}, {\cal E}]$
is divisible on both sides.
If $X\in [\Delta^{o}, {\cal E}]$ and $A\in  {\bf S} $,
then
$$A\backslash X=\int_{[n]\in \Delta} X_n^{A_n}.$$
If  $B\in  {\cal E} $, then $(X/B)_n={\cal E}(B,X_n)$
for every $n\geq 0$.

\medskip

To a map  $u:A\to B$ in ${\bf S}$ and a map $v:S\to T$ 
in ${\cal E}$,
we can associate the map
$$u  \Box'  v:A \Box T  \sqcup_{A \Box S}  B\Box S  \longrightarrow B \Box  T
$$
in $[\Delta^{o}, {\cal E}].$
If $f:X\to Y$ is a map
in $[\Delta^{o}, {\cal E}]$
we then have a map
$$
\langle u  \backslash f\rangle: B \backslash X \longrightarrow B \backslash Y
\times_{A \backslash Y}A \backslash X
$$
in ${\cal E}$ and a map 
$$
\langle f/v \rangle: X/T \longrightarrow Y/T
\times_{Y/S}X/S
$$
in ${\bf S}$.

\medskip

In Reedy theory \cite{Ho}, the object 
$\partial \Delta[n] \backslash X$ is the {\it matching space} $M_nX$.
If $\delta_n$ denotes the inclusion $\partial \Delta[n]\subset \Delta[n]$,
then $\delta_n\backslash X$ is the canonical map
$X_n\to M_nX$.  If $f:X\to Y$ is a map
in $[\Delta^{o}, {\cal E}]$
then $ \langle \delta_n\backslash f\rangle$
is the {\it matching map} 
$$
 X_n\longrightarrow Y_n
\times_{M_nY}M_nX
$$
obtained from the square 
$$
\xymatrix{
X_n \ar[d] \ar[r]& M_nX  \ar[d]\\
Y_n \ar[r]    & M_nY . }
$$

\bigskip

These constructions can be dualised.
Recall that a {\it cosimplicial set} is a covariant functor $\Delta\to {\bf Set}$.
If $A$ is a cosimplicial set and 
$X\in [\Delta^{o}, {\cal E}]$
we shall put
$$A\otimes X=\int^{[n]\in \Delta} A_n\times X_n.$$
Let us denote by $\Delta^c[n]$ the cosimplicial object
$\Delta([n],-):\Delta\to {\bf Set}$.
Its boundary is defined
to be the maximal proper sub-object  $\partial \Delta^c[n]\subset \Delta^c[n]$.
In Reedy theory \cite{Ho}, the object 
$ \partial \Delta^c[n] \otimes X$ is the {\it latching space} $L_nX$.
If $\delta^c_n$ denotes the inclusion $\partial \Delta^c[n]\subset \Delta^c[n]$,
then $\delta_n^c\otimes X$ is the canonical map
$L_nX \to X_n$.  If $f:X\to Y$ is a map
in $[\Delta^{o}, {\cal E}]$,
then the map $ \delta^c_n\otimes' f$
is the {\it latching map}
$$
L_nY
\sqcup_{L_nX} X_n
\longrightarrow Y_n
$$
obtained from the square 
$$
\xymatrix{
L_nX \ar[d] \ar[r]& X_n  \ar[d]\\
L_nY \ar[r]    &  Y_n . }
$$

\medskip

Let $({\cal A},{\cal B}) $ be a weak factorisation system in a finitely
bicomplete category ${\cal E}$.
We shall say that a map $f:X\to Y$ in $[\Delta^{o}, {\cal E}]$
is a {\it Reedy} ${\cal A}$-{\it cofibration} if the latching map $ \delta^c_n\otimes' f$
belongs to ${\cal A}$ for every $n\geq 0$.
We shall say that $f$
is a {\it Reedy}  ${\cal B}$-{\it fibration} if the matching map $\langle \delta_n \backslash f\rangle$
belongs to ${\cal B}$ for every $n\geq 0$.

\begin{theorem}  \label{Reedyfact} Let $({\cal A},{\cal B}) $ be a weak 
factorisation system in a finitely
bicomplete category ${\cal E}$.
If ${\cal A}'$ is the class of Reedy ${\cal A}$-cofibrations in $[\Delta^{o}, {\cal E}]$
and ${\cal B}'$ is the class of Reedy ${\cal B}$-fibrations,
then the pair $({\cal A}',{\cal B}') $ 
is a weak factorisation system.
\end{theorem}

\begin{proposition}\label{Reedyfib1} Let $({\cal A},{\cal B}) $ be a weak 
factorisation system in a 
bicomplete category ${\cal E}$.
Then the following conditions on a map $f\in [\Delta^o,{\cal E}]$ 
are equivalent:
\begin{itemize}
\item $f$ is a Reedy ${\cal B}$-fibration;
\item  the map $\langle u\backslash f\rangle$ belongs to $ {\cal B}$
for every monomorphism $u\in {\bf S}$; 
\item  the map $\langle f/v \rangle$ is a trivial fibration
for every $v\in \cal A$.
\end{itemize}
\end{proposition}

\noindent {\bf Proof}: The implication (ii)$\Rightarrow$(i) is obvious.  
Let us prove the implication (i)$\Rightarrow$(iii).
By \ref{saturatedmonosset} it suffices to show that
we have $\delta_m\pitchfork \langle f/v\rangle$ for every $m\geq 0$.
But the condition $\delta_m\pitchfork \langle f/v\rangle$
is equivalent to the condition $v\pitchfork  \langle \delta_m\backslash f\rangle$ 
by \ref{lift2}.
We have $ \langle \delta_m\backslash f\rangle\in {\cal B}$
since $f$ is a Reedy $ {\cal B}$-fibration by assumption.
Thus, $v\pitchfork  \langle \delta_m\backslash f\rangle$ 
since $v\in \cal A$.
This proves that $ \langle  f/v\rangle$ is a trivial fibration.
Let us prove the implication (iii)$\Rightarrow$(ii).
It suffices to show that we have 
$v\pitchfork \langle u\backslash f\rangle$
for every $v\in \cal A$.
But the condition $v\pitchfork  \langle u\backslash f\rangle$  is equivalent to the condition $u \pitchfork \langle f/v\rangle$
by \ref{lift2}.
This proves the result since $\langle f/v\rangle$
is a trivial fibration by hypothesis.
\QED

\medskip

\begin{corollary}  \label{Reedyfact2} Let $({\cal A},{\cal B}) $ be a weak 
factorisation system in a finitely
bicomplete category ${\cal E}$.
If a map $f:X\to Y$ is a Reedy ${\cal B}$-fibration, then the map $f_n:X_n\to Y_n$
belongs to ${\cal B}$ for every $n\geq 0$
\end{corollary}

\noindent {\bf Proof}: For any simplicial set $A$, we have
$A \backslash f=\langle i_A\backslash f\rangle$, 
where $i_A$ denotes the inclusion $\emptyset \subseteq A$.
Thus, $A \backslash f\in {\cal B}$ by \ref{Reedyfib1}.
This proves the result if we take $A=\Delta[n]$.
\QED

\medskip

Let $\cal E$ be a model category with model structure $({\cal C},{\cal W},{\cal F})$.
We shall say that a map $f:X\to Y$ in  $ [\Delta^{o}, {\cal E}]$
is  {\it level-wise acyclic}
 if the map $f_n:X_n\to Y_n$ 
is acyclic for every $n \geq 0$.
We shall say that $f$ is a {\it Reedy cofibration}
if it is a Reedy ${\cal C}$-cofibration.
We shall say that $f$ is a {\it Reedy fibration}
if it is a Reedy ${\cal F}$-fibration.

\begin{theorem} [Reedy, see \cite{Ree}] \label{Reedymodel1}  Let $({\cal C}, {\cal W},{\cal F})$
be a model structure on a category ${\cal E}$.
Then the category  $ [\Delta^{o}, {\cal E}]$ admits a model structure
$({\cal C}', {\cal W}',{\cal F}')$ in which  ${\cal C}'$
is the class of Reedy cofibrations,
${\cal W}'$ is the class of level-wise acyclic maps
and ${\cal F}'$ is the class of Reedy fibrations.
A map $u:A\to B$ in $ [\Delta^{o}, {\cal E}]$
is an acyclic Reedy cofibration iff the latching map
$\delta_n^c \otimes' u $
is an acyclic cofibration for every $n\geq 0$. 
Dually, a map $f:X\to Y$ 
is an acyclic Reedy fibration iff the map
$\langle \delta_n\backslash f \rangle$
is an acyclic fibration for every $n\geq 0$. 
\end{theorem}

We call $({\cal C}', {\cal W}',{\cal F}')$ the {\it Reedy model structure}
associated to  $({\cal C}, {\cal W},{\cal F})$.

\medskip

\begin{proposition}  \label{Reedycofibration2} Let $ {\cal E}$
be a bicomplete model category. If $u\in {\bf S}$ is monic and $v\in {\cal E}$
is a cofibration, then $u  \Box'  v$
is a Reedy cofibration which is acyclic if $v$ is acyclic.
\end{proposition}

\noindent {\bf Proof}: If $u$ is monic 
and $v$ is a cofibration, let us show that 
the map $u  \Box'  v$
is a Reedy cofibration.
For this, it suffices to show that
we have $(u  \Box'  v)\pitchfork f$ for every acyclic Reedy fibration $f$.
Let $({\cal C},{\cal W},{\cal F})$ be the model structure on $ {\cal E}$.
The map $f$ is a Reedy ${\cal F}\cap {\cal W}$-fibration by  \ref{Reedymodel1}.
Hence the map $ \langle  f/v\rangle$ 
is a trivial fibration by \ref{Reedyfib1},
since $v \in {\cal C}$ and since the pair $({\cal C},{\cal F}\cap {\cal W})$ is a weak
factorisation system.
 Hence we have $u\pitchfork   \langle  f/v\rangle$, since $u$ is monic.
It follows that we have $(u  \Box'  v)\pitchfork f$ 
by \ref{lift2}.  
Let us now show that $u  \Box'  v$
is acyclic if moreover $v$ is acyclic.
For this, it suffices to show that
we have $(u  \Box'  v)\pitchfork f$ for every Reedy fibration $f$. 
The map $ \langle  f/v\rangle$ 
is a trivial fibration by \ref{Reedyfib1}, since $v \in {\cal C}\cap {\cal W}$
and since  
the pair $({\cal C}\cap {\cal W},{\cal F})$  is a weak
factorisation system. Hence we have $u\pitchfork   \langle  f/v\rangle$, since $u$ is monic.
It follows that we have $(u  \Box'  v)\pitchfork f$ 
by \ref{lift2}.  
\QED

The object $[0]$ is terminal in the category $\Delta$.
Hence the functor $i:1\to \Delta$ defined by putting
$i(1)=[0]$ is right adjoint to the constant functor $p:\Delta \to 1$.
It follows that the functor $i^*: [\Delta^{o}, {\cal E}]\to {\cal E}$
is right adjoint to the functor $p^*: {\cal E}  \to [\Delta^{o}, {\cal E}] $.
If $X\in [\Delta^{o}, {\cal E}]$, then  $i^*(X)=X_0$;
if $A\in {\cal E}$, then $p^*(A)=1\Box A$ is the constant simplicial
object with value $A$.

\begin{proposition}  \label{Quillen pair for Reedy} Let $ {\cal E}$
be a model category. Then the pair of adjoint
functors 
$$p^*: {\cal E}  \leftrightarrow [\Delta^{o}, {\cal E}] :i^*$$ 
is a Quillen pair, where  $[\Delta^{o}, {\cal E}] $
is given the Reedy model structure. 
\end{proposition}

\noindent {\bf Proof}: 
It suffices to show that
the functor $p^*$ is a left Quillen functor by \ref{Quillenpair1}.
This follows from \ref{Reedycofibration2}
since $p^*(A)=1\Box A$.

\bigskip 
ÊÊ

\section{Epilogue}

A  Quillen model category is the primary example
 of an ``abstract homotopy theory" \cite{Q}.
 It permits the construction of
homotopy pullbacks and pushouts, 
of fiber and cofiber sequences, etc.
Dwyer and Kan have proposed using general simplicial categories
to model general homotopy theories, see \cite{DK1} and \cite{DK2}.
Recently, Bergner \cite{B1} has established a model structure
on the category of simplicial categories.
In \cite{B2}, she introduces a new model structure ${Segal\ cat}'$ on
the category of Segal precategories 
and she obtains a chain of Quillen equivalences:
$${Simp.cat} \to {Segal\ cat}' \leftarrow {Segal\ cat} \leftarrow {Comp.\ Segal\ sp}$$
(in this chain, a Quillen equivalence is represented by the right adjoint functor).
See \cite{B4} for a survey.
It follows from her results combined the result proved in the present paper that the model structure for simplicial categories is indirectly Quillen equivalent 
to the model structure for quasi-categories.
We shall see in \cite{J4} that the coherent nerve functor
defined by Cordier \cite{C} and studied in Cordier-Porter \cite{CP}
defines a direct Quillen equivalence 
$${Simp.cat}\to {Quasicat}$$
between the model structure for simplicial categories and the
model structure for quasi-categories.
See Lurie \cite{Lu1} for another proof.
An axiomatic approach to proving all the equivalences above
was proposed by T\"oen in \cite{T2}.

\medskip

There are other important notions
of homotopy theories, for example the $A_\infty$-spaces 
of Stasheff \cite{MSS}. See also Batanin  for
$A_\infty$-categories \cite{Ba}. 
A theory of homotopical categories was developed by
Dwyer, Hirschhorn, Kan and Smith \cite{DHKS}.
Dugger studies universal  homotopy theories in \cite{D}.
The model structure for quasi-categories belongs to a 
class of model structures in presheaf categories studied
by Cisinski \cite{Ci}.
A quite different notion of homotopy theory was  introduced
by Heller \cite{He} based on the idea of hyperdoctrine of Lawvere \cite{La}.
A similar notion called {\it d\'erivateur} was later introduced by Grothendieck 
and studied by Maltsiniotis \cite{M1}. See \cite{M2} for an extension.

\bibliographystyle{amsalpha}

\end{document}